\documentclass[11pt,a4paper]{amsart}
\usepackage[headings]{fullpage}
\usepackage{amssymb,amscd,hyperref}
\usepackage{amsthm}
\usepackage[alphabetic, nobysame]{amsrefs}
\usepackage[all,cmtip]{xy}%\CompileMatrices
\usepackage{color}
\usepackage{xcolor}
\usepackage{tikz-cd} 
\numberwithin{equation}{section}
\newtheorem{thm}{Theorem}[section]
\newtheorem{lem}[thm]{Lemma}
\newtheorem{prop}[thm]{Proposition}
\newtheorem{cor}[thm]{Corollary}

\theoremstyle{definition}
\newtheorem{defn}[thm]{Definition}
\newtheorem{ex}[thm]{Example}

\theoremstyle{remark}
\newtheorem{rem}[thm]{Remark}
\newtheorem*{acknowledgements}{Acknowledgements}

\def\Z{\mathbb{Z}}

\def\THH{\mathsf{THH}}
\def\thr{\mathsf{THR}}

\def\mack{\mathsf{Mack}}
\def\tamb{\mathsf{Tamb}}

\def\coind{\mathsf{coind}}
\def\cgrings{\mathsf{cGrings}}
\def\rot{\mathrm{rot}}

\def\fsets{\mathsf{Sets^f}}

\def\bigsquare{\raisebox{-1mm}{\scalebox{1.5}{$\square$}}}

\def\norm{\mathsf{norm}}
\def\trace{\mathsf{tr}}
\def\res{\mathsf{res}}
\newcommand{\und}[1]{{\underline{#1}}}

\newcommand{\cL}{\mathcal{L}}
\newcommand{\cH}{\mathcal{H}}

\def\ra{\rightarrow}
\def\leq{\leqslant}
\def\geq{\geqslant}

\def\cala{\mathcal{A}}

\newcommand{\id}{{\mathrm{id}}}

\newcommand{\ie}{\emph{i.e.}}

\newcommand{\comm}{\mathrm{Comm}}

\def\ie{\emph{i.e.}}

\def\id{\mathrm{id}}

\begin{document}
\title{Loday constructions of Tambara functors}
\author{Ayelet Lindenstrauss}
\address{Mathematics Department, Indiana University, 831 East Third Street,
  Bloomington, IN 47405, USA}
\email{alindens@iu.edu}
\author{Birgit Richter}

\address{Fachbereich Mathematik der Universit\"at Hamburg,
Bundesstra{\ss}e 55, 20146 Hamburg, Germany}
\email{birgit.richter@uni-hamburg.de}

\author{Foling Zou}
\address{Academy of Mathematics and Systems Science, Chinese Academy
of Sciences,  
No. 55 Zhongguancun East Road, Beijing 100190,  China}
\email{zoufoling@amss.ac.cn}

%\date{\today}
\maketitle

\begin{abstract}

Building on work of Hill, Hoyer and Mazur we propose an equivariant version of 
a Loday construction for $G$-Tambara functors where $G$ is an arbitrary finite
group. For any finite simplicial $G$-set and any $G$-Tambara functor, our
Loday construction is a
simplicial $G$-Tambara functor. We study its properties and examples. For a
circle with rotation action by a finite cyclic group our construction agrees
with the twisted cyclic nerve of Blumberg, Gerhardt, Hill, and Lawson. We also show how the Loday construction for genuine commutative $G$-ring spectra relates to our algebraic one via the $\und{\pi}_0$-functor. We describe Real topological Hochschild homology as such a Loday construction. 
  
\end{abstract}

% \keywords{Tambara functors; relative norm construction; Real topological Hochschild homology; twisted cyclic nerve; Loday construction}
\section{Introduction}

In this paper we generalize the Loday
construction for
commutative rings \cite{pirashvili} to the
equivariant context where the groups involved are finite. For every genuine commutative $G$-ring spectrum $R$,
Brun showed \cite{brunpi0} 
that $\und{\pi}_0(R)$ is a $G$-Tambara functor -- this is roughly a
Mackey functor
with compatible commutative ring structures and multiplicative norms; see
\cite[\S 2]{tambara} for the definition. We will define a Loday construction
for $G$-Tambara functors for any finite group $G$ and investigate some of
its properties.

Our work builds on the Hill-Hopkins notion of a $G$-symmetric
monoidal category \cite{hh}. In equivariant algebra there is an important
difference between commutative
monoids in the category of $G$-Mackey functors -- these are the commutative
$G$-Green functors -- and $G$-commutative monoids -- these are the $G$-Tambara
functors. Mazur, Hill-Mazur and Hoyer
\cite{mazur,hm,hoyer}  prove that for any finite group and any
$G$-Tambara functor $\und{R}$ there is a compatible definition of the tensor
product of a finite $G$-set $X$ with $\und{R}$.

Their construction leads directly to our definition of a Loday construction in
Section \ref{sec:def} because their naturality statements ensure that one can  extend this tensor product to a tensor 
product of a $G$-Tambara functor with a finite simplicial $G$-set.  
 We study basic properties like the behaviour of the Loday
construction on free Tambara functors and on norms in Section \ref{sec:basic} and we confirm that restricting a $G$-Loday construction to $H$ for a subgroup $H < G$ gives the $H$-Loday construction of the $H$-restricted Tambara functor. 
We prove in Section \ref{sec:htpinv} that the Loday construction is homotopy
invariant with respect to $G$-homotopy equivalences. Even for constant
Tambara functors the Loday construction detects interesting features
of the $G$-simplicial set and embedding commutative rings into the
equivariant setting by taking their constant Tambara functors doesn't produce any undesired properties. We discuss these facts in Section
\ref{sec:constant}. We show how the $\und{\pi}_0$-functor relates
Loday constructions for genuine $G$-equivariant commutative ring
spectra to the Loday construction of the corresponding Tambara functors 
in Section \ref{sec:pi0}.

Section \ref{sec:examples} is the heart of the paper where we present several
explicit examples of Loday constructions. We relate
our Loday construction for a circle with rotation action to the twisted 
cyclic nerve of \cite{bghl} and show that at the free level we just obtain a
subdivision of the ordinary Loday construction.  We also identify Loday constructions for unreduced
suspensions where  we either fix the suspension apices or a $C_2$-action
flips them. 

We show that Real topological Hochschild homology 
can be expressed as an equivariant Loday construction for 
flat and well-pointed genuine commutative $C_2$-ring spectra. In this example
the finite simplicial  $C_2$-set is a circle with flip action.

For the classical Loday construction working relative to a base-ring is often crucial for performing calculations. We propose a relative version of the Loday construction for Tambara functors in Section \ref{sec:rel}. 

\begin{acknowledgements} We thank Mike Hill, Mike Mandell and Peter
  May for helpful comments. 
The first named author was supported by NSF grant DMS-2004300 and a grant
from the Simons Foundation (917555, Lindenstrauss). The
second named author thanks Churchill College Cambridge for its hospitality. 
\end{acknowledgements}

\section{Definition and basic properties} \label{sec:def}

Let $G$ be a finite group, and let $\und{R}$ be a $G$-Tambara functor, as in \cite[\S 2]{tambara}.  There is a tensor product of $\und{R}$ with finite $G$-sets and
our definition of Loday constructions for $G$-Tambara functors is based on that. In the context of cyclic $p$-groups for a prime $p$, this was developed in Mazur's thesis 
\cite[Theorem 2.3.1]{mazur} and published in joint work of Hill and Mazur
\cite[Theorem 5.2]{hm}. For general finite groups Hoyer provided a construction in \cite[\S 2.4]{hoyer}.

In the following let $G$ be a finite group and we denote by  $\mack_G$
the category of
$G$-Mackey functors. This is a symmetric monoidal category with respect to the
$\Box$-product of Mackey functors. The starting point is a $G$-symmetric monoidal structure
in the sense of
\cite[Definition 3.3.]{hh}. Let $\fsets_G$ denote the category of finite $G$-sets and we consider the wide subcategory where we restrict the morphisms to
isomorphisms of finite $G$-sets. There is a functor 
\begin{equation}
\label{Mackeytensor}
(-)\otimes (-)\colon  (\fsets_G, G\mathrm{-isoms}) \times \mack_G \to \mack_G 
\end{equation}
which satisfies the following properties:
\begin{enumerate}
\item For all $X$ and $Y$ in $\fsets_G$ and  $\und{M}$, $\und{N}$ in $\mack_G$,
  there are natural isomorphisms
\begin{equation}
\label{union}
(X \amalg Y)\otimes \und{M} \cong (X\otimes \und{M}) \Box (Y\otimes \und{M})
\end{equation}
and 
\begin{equation}
\label{box}
X\otimes (\und{M} \Box \und{N}) \cong (X\otimes \und{M}) \Box (X\otimes \und{N}).
\end{equation}

\item There is a natural isomorphism
\begin{equation}
\label{prod}
X\otimes (Y\otimes \und{M}) \cong (X\times Y)\otimes \und{M}.
\end{equation}

\item On the category with objects finite sets with trivial $G$-action and morphisms consisting only of isomorphisms, the functor restricts to exponentiation $X\otimes \und{M}= \bigsquare_{x\in X} \und{M}$.

\end{enumerate}
\medskip

If we now turn to the category of $G$-Tambara functors, $\tamb_G$,
Definition 5.3 and Proposition 5.4 in \cite{hm} for cyclic $p$-groups
and \cite[Theorem 2.7.4]{hoyer} in the general case allow us to extend
the definition in (\ref{Mackeytensor}) to the category $\fsets_G$
where morphisms are all $G$-maps,
\begin{equation} \label{eq:tambaratensor}
(-)\otimes (-)\colon \fsets_G \times \tamb_G \to \tamb_G.\end{equation}

The tensor product from \eqref{eq:tambaratensor} has an explicit description. 
By \cite[Theorem 5.2]{hm} for the case of $G$ a cyclic $p$-group and more generally \cite[Theorem 2.7.4]{hoyer} we get that for every subgroup  $H < G$
tensoring with $G/H$ can be identified as 
\begin{equation}\label{orbittensor}
G/H\otimes \und{R} \cong N_H^G i_H^* \und{R},
\end{equation}
where $i_H^*$ is the restriction functor $i_H^* \colon \tamb_G \ra \tamb_H$ 
that takes a $G$-Tambara functor $\und{R}$ to the $H$-Tambara functor 
\[ i_H^*\und{R}(H/K) = \und{R}(G \times_H H/K) \cong  \und{R}(G/K) \]
for any subgroup $K\leq H$; $N_H^G$ is the norm, as in \cite[Definition 3.13]{hm} and \cite[Theorem 2.1.1]{hoyer}. 
Using \eqref{union}, we immediately get that if $X=\amalg_{\alpha\in A} G/H_\alpha$ for subgroups $H_\alpha\leq G$, then 
\begin{equation}\label{orbittensors}
X\otimes \und{R} \cong \bigsquare_{\alpha\in A} N_{H_\alpha}^G i_{H_\alpha }^* \und{R}.
\end{equation}

\begin{rem}
  Since Tambara functors are formally just a diagram category with entries that are commutative rings, we could use direct limits of commutative rings to give us direct limits of Tambara functors, which still have the required structure for being Tambara functors.  This allows us to extend the functor from
  \eqref{eq:tambaratensor} to $G$-sets that are not necessarily finite: 
\begin{equation}\label{tensoring}
(-)\otimes (-)\colon \mathsf{Sets}_G \times \tamb_G \to \tamb_G.\end{equation}

\medskip
This generalizes a construction in the nonequivariant setting, that is crucial
for the Loday construction. For any finite set $X$ and any  commutative ring
$R$, the assignment
\[ (X, R) \mapsto X\otimes R = \bigotimes_{x\in X} R \]
is functorial in $X$ and $R$. For an arbitrary (possibly infinite set)
$X$ this can be defined  as the colimit over finite subsets of $X$ of this
construction.  Maps $f\colon X\to Y$ between finite sets
are sent to 

\begin{equation*}
f_*\colon \bigotimes_{x\in X} R \to \bigotimes_{y\in Y} R,\quad\quad  f_*(\bigotimes_{x\in X}r_x)=\bigotimes_{y\in Y} b_y\quad \text{ with }  \quad b_y=\prod_{x\in f^{-1}(y)} r_x. 
\end{equation*}
This is well-defined because  $R$ is assumed to be commutative.

An analogous definition works for simplicial rings and for commutative ring
spectra.

This construction is the basis of the (nonequivariant) Loday construction, which takes a simplicial set $X_\cdot$ and a commutative ring  $R$  
and sends them to the simplicial commutative ring

\begin{equation} \label{eq:ordinaryloday}
 \cL_{X_\cdot}(R) =\{[n]\mapsto X_n\otimes R = \bigotimes_{x\in X_n} R \} 
\end{equation}
such that the simplicial structure maps are 
induced by the corresponding maps in $X$.  If $R_\cdot$ is a simplicial commutative ring, then we take the diagonal of the bisimplicial result of the above definition. That results in a simplicial commutative ring. If $R$ is a commutative ring  spectrum, we use the smash product instead of the tensor product and take the
realization of the resulting simplicial commutative ring spectrum to obtain a  commutative ring spectrum. 
\end{rem}

\begin{defn} \label{def:equivariantloday}
Let $G$ be a group, $X_\cdot$ be a simplicial $G$-set, and $\und{R}$ be a $G$-Tambara functor.  \emph{The equivariant Loday construction of $\und{R}$ with respect to $X_\cdot$} is defined to be the simplicial $G$-Tambara functor
\begin{equation*}
\cL^G_{X_\cdot}(\und{R}) =\{[n] \mapsto X_n\otimes \und{R} \}, 
\end{equation*}
where $X_n\otimes \und{R}$ is defined by the functor in \eqref{tensoring}.  The simplicial structure maps are 
induced by the corresponding maps in $X$.
\end{defn}

\begin{rem}
By construction the equivariant Loday construction is natural in the simplicial $G$-set $X_\cdot$ and in the $G$-Tambara functor $\und{R}$. We will mostly consider finite simplicial $G$-sets $X$. These are simplicial objects in finite $G$-sets. 
\end{rem}

\begin{prop} \label{properties}
The equivariant Loday construction satisfies the following properties:
\begin{enumerate}
\item
For all finite simplicial $G$-sets $X_\cdot$ and $Y_\cdot$ and $G$-Tambara functors $\und{R}$ and $\und{T}$, there are natural isomorphisms of simplicial Tambara functors
\begin{equation*}
\cL^G_{X_\cdot\amalg Y_\cdot}(\und{R})\cong \cL^G_{X_\cdot}(\und{R}) \Box \cL^G_{Y_\cdot}(\und{R})
\end{equation*}
and 
\begin{equation*}
\cL^G_{X_\cdot}(\und{R} \Box \und{T}) \cong \cL^G_{X_\cdot}(\und{R})  \Box \cL^G_{X_\cdot}(\und{T}).
\end{equation*}
\item For all finite simplicial $G$-sets $X_\cdot$ and $Y_\cdot$ and any
  $G$-Tambara functor $\und{R}$, there is a natural isomorphism
  between the diagonal of the bisimplicial Tambara functor
  $\cL^G_{X_\cdot}(\cL^G_{Y_\cdot}(\und{R}))$ and the simplicial Tambara functor
  $\cL^G_{X_\cdot\times Y_\cdot}(\und{R})$. 
\item If $X_\cdot$ and $Y_\cdot$ are finite simplicial $G$-sets containing a common
 simplicial $G$-subset $Z_\cdot$ and $\und{R}$ is  $G$-Tambara functor, 
then 
\begin{equation*}
\cL^G_{X_\cdot\sqcup_{Z_\cdot}  Y_\cdot}(\und{R})\cong \cL^G_{X_\cdot}(\und{R}) \Box_{\cL^G_{Z_\cdot}(\und{R}) } \cL^G_{Y_\cdot}(\und{R}).
\end{equation*}

\end{enumerate}

\end{prop}

\begin{proof}
  These are all proved levelwise, using the isomorphisms of \eqref{union},
  \eqref{box}, and \eqref{prod}, and the fact that for all $n\geq 0$,
  \[(X_n \sqcup_{Z_n} Y_n) \otimes \und{R} \cong (X_n\otimes \und{R}) \Box_{Z_n\otimes \und{R}} (Y_n\otimes \und{R})\]
  also using
  \eqref{union}. Note that $Z_\cdot$ has to contain entire orbits to be a simplicial $G$-subset. 
\end{proof}

\begin{rem} \label{rem:Weyl}
  There is an explicit description of the action of the Weyl group of $H$ in
  $G$,
  $W_G(H)$, on terms of the form $N_H^Gi_H^*\und{R}$ for instance in \cite[Proof of Proposition 5.10]{hm} in the case of cyclic $p$-groups. This combines a cyclic permutation action and a coordinatewise action.  It is also observed in
  \cite[Theorem 5.11]{hm} that $N_H^Gi_H^*\und{R}$ with this Weyl action is
  isomorphic to $N_H^Gi_H^*\und{R}$ with the Weyl action from
  \cite[Fact 5.8]{hm}. So for the Loday construction there is a choice to make and we choose to work with the first Weyl action.
\end{rem}

\section{Basic constructions} \label{sec:basic}
\subsection{Free Tambara functors}
For every $G$-Tambara functor $\und{R}$ the commutative ring $\und{R}(G/e)$
carries a $G$-action that is compatible with the ring structure. We call the
category of such rings together with equivariant ring maps the category of
commutative $G$-rings and we denote it by $\cgrings$. 
\begin{lem}(\cite[\S 2]{brunfree})
  \label{lem:brun}
  There exists a left adjoint $\und{F}(-)\colon \cgrings \ra \tamb_G$ to the
  functor
  \[ \tamb_G \ra \cgrings, \quad \und{R} \mapsto \und{R}(G/e).\]
\end{lem}
\begin{proof}
  Brun first shows that the category of commutative $G$-rings is equivalent to
  the category of $fG$-Tambara functors \cite[Lemma 6 (i)]{brunfree}. Here, an $fG$-Tambara functor just accepts free orbits as input. He then constructs the free $G$-Tambara
functor $\und{F}(R)$ as the left Kan extension of the $fG$-Tambara functor corresponding to $R$  along the inclusion $U_{fG} \hookrightarrow U_G$
\cite[p.~241]{brunfree}. Here, $U_{fG}$ ($U_G$)  denotes the category of bispans based on finite free $G$-sets (all finite $G$-sets). 
\end{proof}

    A concrete formula of $\und{F}(-)$ in the case of $G=C_p$ can be found in \cite[Example 1.4.8]{mazur}.

\begin{rem} \label{rem:notfaithful}

Note that the right adjoint functor $U$ that sends a $G$-Tambara functor $\und{R}$ 
to $\und{R}(G/e)$ is not faithful: for simplicity assume that  $G=C_2$ and
let $\und{R}$ be a $C_2$-Tambara functor. Let $M$ be a free $\und{R}(C_2/C_2)$-module of rank $1$ with generator $m$. Then we can define a new $C_2$-Tambara functor $\und{R} \rtimes M$ as $\und{R} \rtimes M(C_2/e) = \und{R}(C_2/e)$ and
$\und{R} \rtimes M(C_2/C_2) := \und{R}(C_2/C_2) \oplus M$ with the square-zero
multiplication. We keep the norm and transfer maps from $\und{R}$ on $\und{R} \rtimes M$ and define the restriction of $m$ to be zero. Then
\begin{align*}
& \cgrings(U(\und{R} \rtimes M), U(\und{R} \rtimes M)) =
                 \cgrings(\und{R} \rtimes M(C_2/e), \und{R} \rtimes M(C_2/e)) \\
  = &   
      \cgrings(\und{R}(C_2/e), \und{R}(C_2/e)) = \cgrings(U(\und{R}), U(\und{R}))
      \end{align*}
whereas we can define at least two different self-maps of the $C_2$-Tambara functor
$\und{R} \rtimes M$ by taking the identity on the free level and one morphism that
sends $(r,m) \in \und{R}\rtimes M(C_2/C_2)$ to $(r,0)$ and a second one that
sends $(r,m)$ to $(r,m)$.

A consequence is that the counit of the adjunction $\varepsilon \colon \und{F} \circ U \ra \mathrm{Id}$ is \emph{not} an epimorphism (see e.g.~\cite[Proposition 2.4.11]{cats}).  
\end{rem}

In the following we denote by $\res_H^G$ the forgetful functor from $G$-sets to $H$-sets for a subgroup $H \leq G$. We recall the following standard isomorphism. Let $Z$ be a finite
$G$-set. Then 
\begin{equation} \label{eq:giso}
G \times_H \res_H^GZ \cong G/H \times Z
\end{equation}
  where the isomorphism is given by 
 $[g,z] \mapsto (gH,gz)$ with inverse $(gH, z) \mapsto [g,g^{-1}z]$.
This isomorphism transforms the $G$-action on the left factor of $G \times_H \res_H^GZ$ to the diagonal $G$-action on $G/H \times Z$.

\begin{lem} \label{lem:adjoint}
    Let $G$ be a finite group, and let $Y$ be a finite $G$-set. 
    Then the functor
    \[ Y \otimes (-) \colon \tamb_G \ra \tamb_G \]
    is a left adjoint. Its right adjoint sends $\und{T}$ to $\und{T}^Y$ which maps a finite $G$-set $Z$ to  
\[ \und{T}^Y(Z) := \und{T}(Y  \times Z). \]
  
\end{lem}
\begin{proof}
Again, it suffices to prove the claim for orbits.  By definition $G/H \otimes \und{R} = N_H^Gi_H^*\und{R}$ and $N_H^G(-)$ is left adjoint to restriction. Therefore for every $G$-Tambara functor $\und{T}$ 
  \[ \tamb_G(N_H^Gi_H^*\und{R}, \und{T}) \cong \tamb_H(i_H^*\und{R},
    i_H^*\und{T}). \]
  But the restriction functor $i_H^*$ is itself a left adjoint and Strickland calls its right adjoint coinduction
  \cite[Prop 18.3]{strickland}. For an $H$-Tambara functor $\und{S}$ and a finite $G$-set $X$ the
  latter is defined as $\coind_H^G\und{S}(X) := \und{S}(\res_H^ G(X))$.  Therefore
  \[ \tamb_H(i_H^*\und{R},
    i_H^*\und{T}) \cong \tamb_G(\und{R}, \coind_H^Gi_H^*\und{T}). \]
We know that $i_H^{*}\und{T}(-) = \und{T}(G \times_H (-))$ and
  $\coind_H^G\und{T}(-) = \und{T}(\res_H^G(-))$. The isomorphism from \eqref{eq:giso} then finishes the proof.  
\end{proof}

On the level of Mackey functors, coinduction is actually the left and right adjoint of restriction $i_H^*$ (see \cite[p.~1871]{tw}, where coinduction is called induction), but on the level of Tambara functors $N_H^G(-)$ is the left adjoint and $\coind_H^G$ is the right adjoint.

  If $T_\cdot$ is a simplicial commutative $G$ ring, we
  define $\und{F}(T_\cdot)$ as $\und{F}(T_\cdot)_n
:= \und{F}(T_n)$. 

\smallskip
Note that for any commutative $G$-ring $R$ and any finite set $X$, $X\otimes R = \bigotimes_{x\in X} R $ is a commutative $G$-ring as well, with $G$ acting on all copies of $R$ simultaneously.
Similarly, for any finite simplicial set $X_\cdot$, the usual Loday construction $\cL_{X_\cdot}(R)=X_\cdot \otimes R $
with
\[ (X_\cdot \otimes R)_n := \bigotimes_{x\in X_n} R\]
with simplicial structure maps as in \eqref{eq:ordinaryloday} and $G$ acting on all copies of $R$ simultaneously 
 is a simplicial commutative $G$-ring.
 
 \begin{rem} \label{rem:coproduct}
For any finite set $X$ and commutative $G$-ring $R$, the coproduct of copies of $R$ indexed by $X$ in the category of commutative rings $X \otimes R$ is also the corresponding coproduct in the category of commutative $G$-rings.  The crucial point is that the inclusion of any copy of $R$ indexed by $x\in X$ into $X\otimes R$ is $G$-equivariant by the choice of action on $X\otimes R$, and because of that compatibility  factoring $G$-equivariant maps from $X$ copies of $R$ through $X \otimes R$ will automatically give a $G$-equivariant map from $X \otimes R$.
\end{rem}

We identify the equivariant Loday construction of free Tambara functors as follows:
\begin{prop}
  For every finite group $G$, for every commutative $G$-ring $R$ and every finite simplicial $G$-set $X_\cdot$ we have
  \[ \cL_{X_\cdot}^G\und{F}(R) \cong \und{F}(\cL_{X_\cdot}(R)),  \]
  where the Loday construction $\cL_{X_\cdot}(R)$ is defined using the underlying simplicial set $X_\cdot$ where we forget the $G$-action.
\end{prop}

\begin{proof}
  It suffices to check the claim degreewise and because of \eqref{union} it
  suffices to check it on orbits $G/H$.  We write $|G/H|$ for the set
    $G/H$ after forgetting the $G$-action. 
  For every $G$-Tambara functor $\und{T}$
  we have
  \begin{align*}
    \tamb_G(G/H\otimes \und{F}(R), \und{T})
     & \cong \tamb_G(\und{F}(R), \und{T}(G \times_H \res_H^G(-))) \qquad \text{ by Lemma \ref{lem:adjoint}} \\
     & \cong \cgrings(R,\und{T}(G \times_H \res_H^GG/e) ) \qquad \text{ by Lemma \ref{lem:brun}} \\
     & \cong \cgrings(R,\und{T}(G/H \times G/e))  \\
     & \cong \cgrings(R,\und{T}(|G/H| \times G/e)) \\
     & \cong \cgrings(R, \prod_{|G/H|}\und{T}(G/e)) \\
     & \cong \cgrings(|G/H| \otimes R, \und{T}(G/e)) \qquad \text{ by Remark \ref{rem:coproduct}} \\
    & \cong \tamb_G(\und{F}(|G/H| \otimes R), \und{T})  \qquad \text{ by Lemma \ref{lem:brun}}. 
  \end{align*}  

For the first and second unlabelled isomorphism we use the isomorphism of $G$-sets from \eqref{eq:giso} 
$G \times_H \res_H^GG/e \cong G/H \times G/e$.   We then compose this
isomorphism with the automorphism that sends $(gH, \tilde{g}e)$ to
$(\tilde{g}^{-1}gH, \tilde{g}e)$. This reduces the diagonal $G$-action to the $G$-action on the right factor $G/e$. Thus in the last three rows of the equation, we regard $G/H$ as a set rather than as a $G$-set.

For the third unlabelled isomorphism, we use that $\und{T}$ is a Tambara functor, so it turns disjoint unions of finite $G$-sets into products. 
\end{proof}

\subsection{Turning a commutative ring into a Tambara functor via the norm}

Let $R$ be a commutative ring. We can view $R$ as an $e$-Tambara functor, where
$e$ denotes the group with one element. For all finite groups $G$, $N_e^GR$ is
then a $G$-Tambara functor.

An immediate consequence of Proposition \ref{properties} (2)  is the following fact:
\begin{prop}
  Let $R$ be a commutative ring and let $\und{R}$ be any $G$-Tambara functor with
  $i_e^*\und{R} = \und{R}(G/e) =  R$. Then for all finite simplicial $G$-sets
  $X_\cdot$ 
  \[ \cL_{X_\cdot}^G (N_e^GR) \cong \cL_{X_\cdot}^G (N_e^Gi_e^*\und{R})  \cong
    \cL^G_{X_\cdot} (G/e \otimes \und{R}) \cong \cL^G_{X_\cdot \times G}(\und{R}) \]
  where we view $G$ as the constant simplicial set on the right. 
\end{prop}

We also apply $N_e^G$ to a simplicial commutative
ring degreewise and obtain a simplicial $G$-Tambara functor. If we apply $N_e^G$ to the non-equivariant Loday construction of a commutative ring, we get
the following relationship: 

\begin{prop}
  Fix a commutative ring $R$. Let $\und{R}$ be any $G$-Tambara functor with
  $\und{R}(G/e) = R$ and let $X_\cdot$ be any finite simplicial set. Then
  \[ N_e^G\cL_{X_\cdot}(R) \cong \cL_{G \times X_\cdot}^G(\und{R}).   \]
\end{prop}
\begin{proof}
  The claim follows directly from Hoyer's naturality requirement in
  \cite[Definition 2.7.2]{hoyer} that he proves in \cite[Theorem 2.7.4]{hoyer}:
  For any finite group $G$, every pair of
  subgroups $H \leq K \leq G$  and every finite $H$-set $S$ there is an isomorphism of $K$-Tambara functors
  \begin{equation} \label{eq:norminduction}
(K \times_H S) \otimes (i_K^*\und{R}) \cong N_H^K(S \otimes (i_H^*\und{R}))
\end{equation}
which is natural in the $H$-set $S$. This implies that on the level of Loday constructions we obtain for every finite simplicial $H$-set $X_\cdot$
\[ \cL_{K \times_H X_\cdot}^K(i_K^*\und{R}) \cong N_H^K(\cL_{X_\cdot}^H (i_H^*\und{R})). \]
We apply his result in the situation where $H = e$ and $K=G$. In this case, we
identify $\cL_{G \times_e X_\cdot}^G(i_G^*\und{R}) = \cL_{G \times X_\cdot}^G(\und{R})$ with
$N_e^G(\cL_{X_\cdot} (i_e^*\und{R}))$.  As $i_e^*\und{R}(e/e) = \und{R}(G/e) = R$, this proves the claim. 
\end{proof}  
\begin{rem}
  Note that the norm $N_e^G$ of a commutative ring $R$ does \emph{not} agree with the
  free $G$-Tambara-functor, if we start with a commutative $G$-ring $R$ viewed as a commutative $G$-ring with
  trivial $G$-action: By adjunction $\tamb_G(N_e^GR, \und{T}) \cong \tamb_e(R,\und{T}(G/e))$ where we view $\und{T}(G/e)$ just as a commutative ring.
  In contrast, $\tamb_G(\und{F}(R), \und{T}) \cong \cgrings(R, \und{T}(G/e))$;
    so here, we \emph{do} remember the $G$-action on the commutative ring
  $\und{T}(G/e)$.  
  As $R$ carries the trivial action this yields ring morphisms from $R$ into the $G$-fixed points of $\und{T}(G/e)$, $\und{T}(G/e)^G$. 
\end{rem}

\subsection{Change of groups}
Hoyer's naturality result from \eqref{eq:norminduction} immediately gives the
following isomorphism:

\begin{prop}
  For every finite group $G$ with a  subgroup $H < G$, and every finite
  simplicial $H$-set $X_\cdot$ 
  \[ \cL_{G \times_H X_\cdot}^G(\und{R}) \cong N_H^G\left(\cL_{X_\cdot}^H(i_H^*\und{R})\right) \]
as simplicial $G$-Tambara functors. 
\end{prop}

We also have an isomorphism on Loday constructions with respect to
restrictions. This follows from the next
Theorem \ref{thm:restriction}.
  \begin{prop}
  For every finite group $G$ with a  subgroup $H < G$, and every finite
  simplicial $G$-set $X_\cdot$ 
\[ i^{*}_{H}\cL_{X_\cdot}^G(\und{R}) \cong \cL_{i_H^*X_\cdot}^H(i_H^*\und{R})\]
as simplicial $H$-Tambara functors. 
\end{prop}

\begin{thm}
  \label{thm:restriction}
    For every finite $G$-set $S$ and $G$-Tambara functor $\und{R}$ there is isomorphism 
\begin{equation}
\label{eq:restriction}
 i_H^*(S \otimes \und{R}) \cong  i_H^*(S) \otimes i_H^*(\und{R})
\end{equation}
natural in $S$ and $\und{R}$.
  \end{thm}

  \begin{proof}
    In Hoyer's work \cite[Theorem 2.5.1]{hoyer}, it is shown that
    there is an 
    isomorphism
\[i_H^*(G/K \otimes \und{R}) \cong \bigsquare_{\gamma \in H\backslash G/K} N_{H \cap
    {}^\gamma K}^H(i_{H \cap
    {}^\gamma K}^*\und{R}) \] 
    that is natural in $\und{R}$.
 The naturality in the finite $G$-set is not stated because it only
 makes sense when 
 $\und{R}$ is a $G$-commutative monoid and Hoyer proved the isomorphism in
 greater generality for $G$-Mackey functors. However,
 the arguments in his proof
 also show the functoriality in finite $G$-sets for $G$-Tambara functors:

 We first recall some definitions from \cite{hoyer}.
 Let $A_G$ be the Burnside category of isomorphism classes of spans of finite
 $G$-sets. Then $G$-Mackey functors are product preserving functors $A_G \to
 \mathrm{Set}$ such that the image is levelwise grouplike. For $f
 \colon
 X \ra Y$ a map of finite $G$-sets, there are certain maps $\res_f\in
 A_G(Y,X)$ and $\trace_f \in 
 A_G(X, Y)$, which give the restriction and transfer maps for Mackey functors.
 A span $X \overset{f}{\longleftarrow} Z \overset{g}{ \longrightarrow
 } Y$ is equivalent to the composite $\trace_{g}\res_f$. 
 
 The norm functor for Mackey functors $\mack_K \longrightarrow
 \mack_G$  is the left Kan extension along
 $\mathrm{Map}_{K}(G,-) \colon A_K \longrightarrow A_G$. Using the
explicit formula for the left Kan extension via coends, we obtain that for
$\und{M} \in \mack_K$ and for a finite $G$-set $Z$,  $N^G_K
\und{M}(Z)$ can be expressed as 
\begin{equation*}
N^G_K \und{M} (Z) \cong \int^{X \in A_K} A_G(\mathrm{Map}_K(G,X), Z)
\times \und{M}(X). 
\end{equation*}
We can therefore represent elements of  $N^G_K
\und{M}(Z)$ as 
 \begin{equation}
   \label{eq:coend-1}
\bigg(\mathrm{Map}_K(G,X) \overset{f}{ \longleftarrow } W \longrightarrow Z, x \in \und{M}(X)\bigg)
\end{equation}
up to coend identification, where $W$ is a finite $G$-set and $X$ is a
finite $K$-set. Note that via the adjunction $(i_K,
\mathrm{Map}_K(G,-))$, the $G$-map $f$ factors as
\[ \xymatrix@1{ W \ar[r]^(0.3){\eta} & \mathrm{Map}_K(G,i^*_KW) \ar[rrr]^{\mathrm{Map}_K(G,\bar{f})} & & & \mathrm{Map}_K(G,
X)}, \]
where in the second map $\bar{f}
\colon 
i^*_KW \longrightarrow X$ is the adjoint of $f$. Note that for a
finite $G$-set $W$, the restriction $i_K^*W$ just views $W$ as a
$K$-set. 
So in the coend, using the morphism
$\res_{\bar{f}} \in A_{K}$, the element \eqref{eq:coend-1} is
identified with the following element in ``standard form''
\begin{equation*}
\bigg(\mathrm{Map}_K(G,i^*_KW) \overset{\eta}{\longleftarrow} W \longrightarrow Z, w=\res_{\bar{f}}(x) \in \und{M}(i^*_KW)\bigg).
\end{equation*}
Note that the left map corresponding to the restriction map in the
standard form is always the counit $\eta$. 

 To obtain elements in $i^*_H(G/K \otimes \und{R}) = i^*_H N_K^G i^*_K
 \und{R}$, we use
 that the restriction functor $i_H^* \colon
 \mack_G \longrightarrow \mack_H$ is precomposition with
 $G\times_H-$. Then, the standard form of an element in $i^*_H N_K^G
 i^*_K \und{R}(Y)$ for an
 $H$-set $Y$ is
 \begin{equation}
   \label{eq:coend-standard}
\bigg(\mathrm{Map}_K(G,i^*_KW) \longleftarrow W \longrightarrow G \times_H Y, a \in \und{R}(G \times_{K} i^*_KW)\bigg).
\end{equation}
Up to isomorphism, $W \cong G\times_H B$ for some $H$-set $B$, and
 \eqref{eq:coend-standard}  can be rewritten as
 \begin{equation}
  \label{eq:coend-2}
\bigg(\mathrm{Map}_K(G,i^*_K(G\times_HB)) \longleftarrow G\times_HB \overset{G\times_Hf}{ \longrightarrow }G \times_H Y, a \in \und{R}(G \times_{K} i^*_K(G\times_HB))\bigg).
\end{equation}

To map this element to $i^*_H(G/K) \otimes i^*_H \und{R} \cong \bigsquare_{\gamma \in H \backslash G/K} N^H_{H
  \cap^{\gamma}K} i^*_{H \cap^{\gamma}K}\und{R}$, Hoyer shows that
there is an isomorphism of $K$-sets
\begin{equation*}
\hat{\theta}\colon \coprod_{\gamma} {\gamma^{-1}} \cdot {}^{\gamma}K
\times_{H \cap^{\gamma}K} i^*_{H \cap^{\gamma}K} B \cong i^*_{K} (G \times_H  B).
\end{equation*}
This produces an isomorphism of $G$-sets
\begin{equation}
  \label{eq:coend-G-sets}
   \bar{\theta} = G \times_K \hat{\theta}\colon  \coprod_{\gamma}G \times_{H \cap^{\gamma} K} i^*_{H \cap^{\gamma}K} B \cong G \times_{K} i^*_K(G \times_H B).
 \end{equation}
 As the left map in \eqref{eq:coend-2} does not carry any extra
 information, such elements are in bijection with 
\begin{equation}
\label{eq:coend-3}
\bigg(\prod\mathrm{Map}_{H \cap^{\gamma} K}(H,i^*_{H \cap^{\gamma} K}B) \longleftarrow B \overset{f}{ \longrightarrow } Y, \res_{\bar{\theta}}(a) \in
\prod\und{R}(G \times_{H \cap^{\gamma} K} i^*_{H \cap^{\gamma}K} B)\bigg).
\end{equation}
This is the standard form of elements in $ \bigsquare_{\gamma} N^H_{H \cap^{\gamma}K} i^*_{H \cap^{\gamma}K}\und{R} (Y) $.
After checking that this is well-defined for choices of standard forms,
Hoyer proved the isomorphism \eqref{eq:restriction}.

To prove the functoriality, we first consider the case where $K < 
L$ are subgroups of $G$ and $G/K \ra G/L$ is the quotient map. The isomorphism $\bar{\theta}$ is
compatible with restrictions in the sense that the following diagram
commutes, whose vertical maps are quotients:
\begin{equation*}
  \begin{tikzcd}
    \coprod_{\gamma \in H \backslash G/K}G \times_{H \cap^{\gamma} K} i^*_{H \cap^{\gamma}K} B  \ar[r, "\bar{\theta}"] \ar[d]
    & G \times_{K} i^*_K (G \times_H B) \ar[d]\\
   \coprod_{\gamma' \in H \backslash G/L}G \times_{H \cap^{\gamma'} L} i^*_{H \cap^{\gamma'}L} B \ar[r, "\bar{\theta}"'] & G\times_{L} i^*_L (G \times_H B).
  \end{tikzcd}
\end{equation*}
This proves that for the projection $G/K \ra G/L$, the induced map
\[ i^*_H(G/K) \otimes i^*_H(\und{R}) \longrightarrow i^*_H(G/L)
  \otimes i^*_H(\und{R})\]  via the isomorphism
\eqref{eq:restriction} 
is induced by the projection
\[ \coprod_{\gamma \in H \backslash G/K}G \times_{H \cap^{\gamma} K} i^*_{H \cap^{\gamma} K}
  B \longrightarrow \coprod_{\gamma' \in H \backslash G/L}G \times_{H
    \cap^{\gamma'} L} i^*_{H \cap^{\gamma'}L} B. \]
In other words, it is induced by tensoring $i^*_{H}(\und{R})$ with the
restriction $i^*_H$ of $G/K \longrightarrow G/L$. The functoriality with respect to conjugation
$G/K \longrightarrow G/{}^{\gamma}K$, the initial morphisms 
$\varnothing \longrightarrow G/K$ and the fold maps $\nabla \colon G/K \coprod G/K \longrightarrow G/K$ follows similarly.
\end{proof}

\section{Homotopy invariance} \label{sec:htpinv}

For a cofibrant commutative ring spectrum $A$, the non-equivariant
Loday construction $\cL_{X_\cdot}(A)$ is, by \cite[Chapter
VII.3]{ekmm}, a homotopy invariant of $|X_\cdot|$, and therefore, so
is $\pi_*(\cL_{X_\cdot}(A))$. We prove the $G$-homotopy invariance for
Loday constructions of $G$-spectra later in Proposition
\ref{prop:spectra-homotopy-inv}.

It follows from \cite[Theorem 2.4]{pirashvili} that for two finite
simplicial sets $X_\cdot, Y_\cdot$ whose homology groups are
isomorphic as graded cocommutative $k$-coalgebras, the algebraic Loday
constructions $\cL_{X_\cdot}^k(A)$ and $\cL_{Y_\cdot}^k(A)$ have
isomorphic homotopy groups 
if $A$ is a commutative $k$-algebra and $k$ is a field. We expect an
analogous $G$-homotopy invariance result for
the equivariant version $\cL^G_{X_.} (\und{R})$ where $X_\cdot$ is  a finite
simplicial  $G$-set and $\und{R}$ a $G$-Tambara functor. For now, we only prove
a weaker result:

\begin{thm} \label{thm:htpinvariance}
  Let $X_\cdot$ and $Y_\cdot$ be two finite simplicial $G$-sets, equipped with simplicial
  $G$-maps $f \colon X_\cdot\to Y_\cdot$ and $g \colon  Y_\cdot\to X_\cdot$ and simplicial
  $G$-homotopies $f\circ g \simeq \id_{X_\cdot}$ and $g\circ f \simeq\id_{Y_\cdot}$, where $G$
  acts trivially on $\Delta^1 = \Delta(-,[1])$. 
  Then for any $G$-Tambara functor $\und{R}$, there is a homotopy equivalence
\[\cL^G_{X_.} (\und{R}) \simeq \cL^G_{Y_.} (\und{R}).  \]
\end{thm}

This is, in the usual way, a corollary of the following result. 
\begin{prop} \label{prop:homotopy-invariance-simp}
  Let $f, g\colon  X_\cdot\to Y_\cdot$ be two simplicial $G$-maps between two finite simplicial  $G$-sets, and assume that there is a simplicial $G$-homotopy $\cH \colon X_\cdot\times \Delta^1\to Y_\cdot$  between them, where $G$ acts trivially on $\Delta^1$,
  with $\cH(x, s_0^n(0)) = f(x)$ and $\cH(x, s_0^n(1)) = g(x)$ for all $n\geq 0$, $x\in X_n$.  Then there is a homotopy between $f_*, g_*\colon \cL^G_{X_.} (\und{R})\to \cL^G_{Y_.} (\und{R})$.  
\end{prop}

A proof in the non-equivariant context can for instance be found in \cite[p.~3]{anderson}. 
\begin{proof}
  We view $\cL_{(-)}^G(\und{R})$ as a functor from finite simplicial $G$-sets to
  simplicial $G$-Tambara functors. 
  We use an assembly map. For $t_n \in \Delta([n],[1])$  we get a map $X_n \ra
  X_n \times \{t_n\}$ sending $x_n$ to $(x_n,t_n)$ and hence
  \[ \cL_{X_n}^G(\und{R})
    \ra\cL^G_{X_n \times \{t_n\}}(\und{R}). \]
  This assembles into a $G$-equivariant simplicial map 
  \[\cala \colon \cL^G_{X_.} (\und{R}) \times\Delta(-,[1]) \to
    \cL^G_{X_.\times \Delta^1} (\und{R}). \]
Here, $G$ acts trivially on $\Delta(-,[1])$. 
Finally, we compose 
\[ \cL^G_{\cH}(\und{R}) \circ \cala \colon   \cL^G_{X_.} (\und{R})  \times\Delta^1 \ra  \cL^G_{Y_.} (\und{R})\]
to get the desired homotopy.
\end{proof}

\begin{rem}
  Note, that there is no simplicial $G$-Tambara structure on
  $\cL^G_{X_.} (\und{R}) \times\Delta(-,[1])$, so the homotopy as a map 
  $\cL^G_{X_.} (\und{R})  \times\Delta^1 \ra  \cL^G_{Y_.} (\und{R})$ cannot be a
  morphism in this category. However, we started with maps $f,g$ of simplicial
  $G$-sets, so they induce morphisms of simplicial $G$-Tambara functors. In
  the situation of Theorem \ref{thm:htpinvariance} they then induce an isomorphism on homotopy groups.

\end{rem}  

\begin{ex}
  If $X_\cdot$ is any simplicial $G$-set and $\und{R}$ is a $G$-Tambara functor, then
  there is a homotopy equivalence $\cL^G_{CX_.}(\und{R}) \simeq \und{R}$ where $CX_\cdot$ is the cone on $X_\cdot$ with $G$-action induced from that on $X_\cdot$ (fixing the cone point) and $\und{R}$ is interpreted as the constant simplicial object at $\und{R}$. This holds because $CX_\cdot$ is $G$-simplicially homotopy equivalent to a
  point.

\end{ex}
\section{Calculations with constant Tambara functors} \label{sec:constant}

\begin{defn}
Let $G$ be a finite group.
For any commutative ring $R$ we denote by $\und{R}^c$ the constant $G$-Tambara functor with $\und{R}^c(G/H) = R$ for any subgroup $H\leq G$, with $\norm\colon  \und{R}^c(G/H) \to \und{R}^c(G/K)$ given by  $\norm(a) = a^{[K:H]}$ and  $\trace\colon  \und{R}^c(G/H) \to \und{R}^c(G/K)$ given by  $\trace(a) = {[K:H]}\cdot a$ for all $H \leq  K$ and all $a \in R$, and all restriction maps equal to the identity. The action of all Weyl groups in $\und{R}^c$ is trivial. We can similarly define the constant $G$-Tambara functor on a simplicial commutative ring $R_\cdot$ as $(\und{R_\cdot})^c_n := \und{R_n^c}$. 
\end{defn}

The following result is a sanity check about importing non-equivariant objects
into the equivariant setting: 
\begin{prop} \label{prop:trivial-constant}
  Let $X_\cdot$ be a simplicial set with trivial $G$-action and let $R$ be any
  commutative ring. Then 
\[\cL^{G}_{X_\cdot}(\und{R}^c) \cong \und{\cL_{X_\cdot}(R)}^c, \]
where $\cL_{X_\cdot}(R)$ is the nonequivariant Loday construction from \eqref{eq:ordinaryloday}.
\end{prop}
\begin{proof}
  By  \cite[Lemma 5.1]{lrz} we know that the box product of two constant Tambara functors corresponding to commutative rings is the constant Tambara functor corresponding to the tensor product of these rings. So in every simplicial degree
  $n$ of 
$\cL^{G}_{X_\cdot}(\und{R}^c)$, we have $\Box_{x\in X_n} \und{R}^c \cong \und{(\otimes_{x\in X_n} R)}^c$.
\end{proof}

\medskip

We now restrict our attention to cyclic groups $G=C_p$ for $p$ a prime
and  consider the Burnside Tambara functor which is the initial object in the category $\tamb_{C_p}$.  The Burnside Tambara functor of the trivial group $\{e\}$ is just $\und{\Z}^c$ and this in
turn can be identified with the commutative ring $\Z$. As $\Z$ is the initial
object in the category of commutative rings and as the norm functor is a left adjoint, it sends initial objects to initial objects, thus  
\begin{equation}\label{Aasnorm}
N_e^{C_p}i_e^*(\und{\Z}^c) \cong N_e^{C_p}\Z \cong \und{A},
\end{equation}
where $\und{A}$ denotes the $C_p$-Burnside Tambara functor. Note that $\und{A}$ is  the unit with respect to the box product of $C_p$-Mackey 
functors, and that $i_e^*$ of the constant $C_p$-Tambara functor $\und{\Z}^c$ is the constant $\{e\}$-Tambara functor $\und{\Z}^c$ which is just the commutative ring $\Z$.

In the following example, the geometric properties of a finite
simplicial $C_p$-space determine the behaviour of the Loday
construction. The existence of fixed points decides about the Loday
construction:  
    \begin{prop} \label{Zconst}
If $X_\cdot$ is a finite simplicial $C_p$-set, then  
     \[\cL^{C_p}_{X_\cdot}(\und{\Z}^c) \cong \begin{cases}  \und{\Z}^c, & \text{ if } X_\cdot^{C_p} \neq \emptyset, \\
    \und A,  & \text{ if } X_\cdot^{C_p} =\emptyset, \end{cases}\]
    where in both cases $\und{\Z}^c$ and $\und{A}$ denote simplicial Tambara functors that are constant as simplicial objects.
    \end{prop}

    \begin{proof}
If $X_\cdot^{C_p} = \emptyset$, then at every simplicial level $X_n$ is a 
finite disjoint union of free orbits, $X_n = \bigsqcup_{e \in E_n} C_p/e$,  and therefore 
\[ X_n\otimes \und{\Z}^c\cong \bigsquare_{x \in E_n} N_e^{C_p}i_e^*\und{\Z}^c. \]
By \eqref{Aasnorm} each factor is the Burnside Tambara functor for $C_p$, and
as this is the unit for the box product, we obtain
\[ X_n \otimes \und{\Z}^c\cong \und{A}. \]
As the simplicial structure maps induce multiplication and insertion of units, we obtain that in this case the equivariant Loday construction is isomorphic to
the constant simplicial Tambara functor with value $\und{A}$. 
 
 \smallskip
 If $X_\cdot^{C_p} \neq \emptyset$, there is an $n\geq
 0$ and some $x\in X_n$ which is fixed under $C_p$. Applying iterations of $d_0$ and $s_0$ then yields a fixed point in every simplicial degree. 
 Therefore for all $n\geq 0$,
 \[ X_n \cong \bigsqcup_{x \in E_1} C_p/C_p \sqcup \bigsqcup_{x \in E_2} C_p/e\]
 with $|E_1| \geq 1$ and we obtain 
 \[ X_n\otimes \und{\Z}^c\cong (\bigsquare_{x \in E_1} \und{\Z}^c) \Box
   (\bigsquare_{x \in E_2} \und{A}).  \]
 As $\und{A}$ is the unit for the box product and as $\und{\Z}^c \Box \und{\Z}^c \cong \und{\Z}^c$ by \cite[Lemma 5.1]{lrz} we obtain
 \[ X_n \otimes \und{\Z}^c \cong \und{\Z}^c,\]
 because there is at least one trivial orbit.  Again, the simplicial structure
 maps  induce the identity under this isomorphism.
 \end{proof}

 Note that $\und{\Z}^c$ is \emph{not} the unit for the box product of Mackey
 functors. Despite this fact, we obtain the following result:

 \begin{cor}\label{Zalg}
    For any commutative ring $R$ and any simplicial $C_p$-set $X_\cdot$ for which $X_\cdot^{C_p} \neq\emptyset$,
    \[\cL^{C_p}_{X_\cdot}(\und{R}^c) \cong \und{\Z}^c \Box \cL^{C_p}_{X_\cdot}(\und{R}^c) .\]
    \end{cor}

\begin{proof}
  Since $R\cong\Z\otimes R$, by  \cite[Lemma 5.1]{lrz} we get that $\und{R}^c \cong \und{\Z}^c \Box \und{R}^c$. The claim follows by  applying part (1) of
  Proposition \ref{properties}   and Proposition \ref{Zconst}.
    \end{proof}

\section{The linearization map } \label{sec:pi0}
Several authors have observed, that one can form Loday constructions
for $G$-equivariant commutative ring spectra. One approach uses the
fact that the category of $G$-spectra can be
turned into a $G$-symmetric monoidal structure using the
Hill-Hopkins-Ravenel norm construction \cite[Example
2.1.1]{mazur} such that $G$-equivariant commutative ring
spectra are the $G$-commutative monoids \cite[Corollary
17.4.35]{hill-overview}. A different approach \cite{bds} uses the fact
that for orthogonal spectra the category of objects with $G$-action is
equivalent to genuine $G$-equivariant spectra. This makes it possible to use the
classical Loday construction for commutative ring spectra and to endow
it with $G$-actions. 

If we work with the first approach, then we can form a Loday
construction $\cL^G_X(R)$, of a 
$G$-equivariant commutative ring spectrum $R$ with respect to a finite
simplicial 
$G$-set $X$. This is based on tensoring such ring spectra with finite
$G$-sets,  such that on orbits $G/H$ we obtain
\[ G/H \otimes R \cong N_H^Gi_H^* R. \]
Here, $i_H^*$ denotes the spectral version of the restriction functor
and $N_H^G$ is the HHR-norm \cite[Definition A.52]{hhr}. 

\medskip
We first prove the $G$-homotopy invariance for Loday constructions of
$G$-spectra:  
\begin{prop} \label{prop:spectra-homotopy-inv}
  Let $f\colon  X_\cdot\to Y_\cdot$ be a morphism of 
  simplicial $G$-sets, such that the realization
  $|f| \colon |X_{\cdot}| \to |Y_{\cdot}|$ is a homotopy equivalence,
  and let $R$ be a $G$-equivariant commutative ring spectrum.
      Then $f_* \colon \cL^G_{X_.} (R)\to \cL^G_{Y_.} (R)$ induces an
      equivalence on geometric realizations.  
  \end{prop} 
  \begin{proof}
Following \cite[Sec 2.3.1]{hhr}, we use $\comm^G$ to denote the
category of $G$-commutative ring spectra with $G$-commutative maps,
and in contrast use $\comm_G$ to denote the category with the same
objects but with non-equivariant multiplicative maps.
Then $\comm^G$ is enriched in  spaces and $\comm_G$ is enriched in 
$G$-spaces. 
For two  arbitrary  $G$-equivariant
commutative ring spectra $R$ and 
$T$ and for any finite $G$-set $X$ there is a homeomorphism \cite[\S
2.3.1, p.~25]{hhr}
\begin{align} \label{eq:1}
  \comm^G(X \otimes R, T)& \cong G\text{Top}(X, \comm_G(R, T)). 
\end{align}
Furthermore, when $Z$ is a $G$-space with trivial action, there is a
chain of $G$-equivariant homeomorphisms 
\begin{equation}
\label{eq:4}
 \underline{\text{Top}_G}(Z, \comm_G(R, T)) \cong \comm_G(Z \otimes R, T)
 \cong \comm_G(R, T^Z).
\end{equation}
Here, $\underline{\text{Top}_G}$ means the $G$-space of non-equivariant maps.

With these preparations in place we prove the claim. Let
$\mathbf{S}_{\bullet}: \comm_G \to s\comm_G$ be the total singular
complex functor that sends $T$ to $(\mathbf{S}_{\bullet} T)_n =
T^{\Delta^n}$. 
By adjunction we get the following chain of homeomorphisms 
\begin{align*}
\comm^{G}(|\cL^G_{X_.} (R)|, T)& \cong s\comm^{G}(\cL^G_{X_.} (R), \mathbf{S}_{\bullet}T) \\
      & \cong sG\text{Top}(X_{\cdot}, \comm_G(R, \mathbf{S}_{\bullet}T)) & \text{ by \eqref{eq:1}} \\
      & \cong sG\text{Top}(X_{\cdot},  \mathbf{S}_{\bullet}\comm_G(R, T))& \text{ by \eqref{eq:4}} \\
      & \cong G\text{Top}(|X_{\cdot}|, \comm_G(R, T))
\end{align*}
The Yoneda lemma then implies the claim. 
  \end{proof}

We can compare $\und{\pi}_0$ of the spectral Loday construction for
$R$ to the
Loday construction on the $G$-Tambara functor $\und{\pi}_0(R)$:

\begin{prop} \label{prop:pi0}
Let $X$ be a finite simplicial $G$-set and let $R$ be a connective 
$G$-equivariant commutative ring spectrum. Then there is an
isomorphism of simplicial $G$-Tambara functors:
\[ \und{\pi}_0(\cL_X^G(R)) \cong \cL_X^G(\und{\pi}_0(R)). \]
\end{prop}

\begin{proof}
  We consider the $n$th simplicial degree $\und{\pi}_0(\cL_X^G(R))_n = \und{\pi}_0(X_n \otimes R)$
  and we decompose $X_n$ into $G$-orbits: $X_n \cong G/H_1 \sqcup \ldots
  \sqcup G/H_k$. Then, as $R$ is a $G$-commutative monoid,
  \[ X_n \otimes R \cong (G/H_1 \otimes R) \wedge \ldots \wedge (G/H_k
    \otimes R) \cong (N^G_{H_1}i_{H_1}^*(R)) \wedge \ldots \wedge
    (N^G_{H_k} i_{H_k}^*(R)). \]
  Hoyer \cite[\S 2.3.2]{hoyer} proves that there are natural
  isomorphisms $\und{\pi}_0N_H^GH(\und{M}) \cong
  N_H^G\und{M}$ for any $G$-Mackey functor $\und{M}$, where  the
  latter is the norm functor that he constructed for Mackey functors and $H$ denotes the equivariant Eilenberg-MacLane spectrum. 

If $E$ is a connective $G$-equivariant spectrum, then by \cite[Lemma 5.11]{ullman} the map $E \ra
H\und{\pi}_0(E)$ induces an isomorphism
\[ \und{\pi}_0(N_H^GE) \cong \und{\pi}_0(N_H^G H(\und{\pi}_0E)). \]

  As we assume $R$ to be connective and as the 
  smash factors $N^G_{H_j}i_{H_j}^*(R)$ are connective as well, we
  obtain
  \begin{align*}
\und{\pi}_0((N^G_{H_1}i_{H_1}^*(R)) \wedge \ldots \wedge
    (N^G_{H_k} i_{H_k}^*(R))) &\cong \und{\pi}_0(N^G_{H_1}i_{H_1}^*(R))
    \Box \ldots \Box \und{\pi}_0(N^G_{H_k}i_{H_k}^*(R)) \\
& \cong N^G_{H_1}i_{H_1}^*\und{\pi}_0(R) \Box \ldots \Box
                              N^G_{H_k}i_{H_k}^*\und{\pi}_0(R). 
  \end{align*}  
  Therefore, in every simplicial degree, 
\[ \und{\pi}_0(\cL_X^G(R))_n \cong (\cL_X^G(\und{\pi}_0(R)))_n. \] 

The simplicial structure maps induce morphisms that come from the
norm-restriction  adjunction or that are induced by the multiplicative
structure on $R$ and $\und{\pi}_0(R)$. As $\und{\pi}_0$ is strong
symmetric monoidal for connective spectra and as the norm
and restriction functors are also strong symmetric monoidal, we get an
isomorphism of simplicial $G$-Tambara functors. 
\end{proof}

Ullman showed that there is no lax symmetric Eilenberg-MacLane
spectrum functor from the category of $G$-Mackey functors to the
category of connective $G$-spectra, if $G$ is a non-trivial finite
group \cite[Theorem 6.2]{ullman}. However, for a $G$-Tambara functor
$\und{R}$ the multiplication on $H\und{R}$ is still induced by the
multiplication map $\und{R} \Box \und{R} \ra \und{R}$:  Ullman
shows \cite[Theorem 5.2]{ullman} that for a connective commutative
$G$-ring spectrum $E$ and a $G$-Tambara functor
$\und{R}$ the morphisms in the homotopy category of commutative
$G$-ring spectra from 
$E$ to $H\und{R}$ are in bijection with the maps of $G$-Tambara
functors from $\und{\pi}_0E$ to $\und{R}$. For $E = H\und{R} \wedge
H\und{R}$, we have $\und{\pi}_0(H\und{R} \wedge H\und{R}) \cong
\und{R} \Box \und{R}$ and therefore the multiplication $\und{R} \Box
\und{R} \ra \und{R}$ gives rise to a multiplication $H\und{R} \wedge
H\und{R} \ra H\und{R}$. In particular, for all $G$-Tambara functors $\und{R}$ Proposition \ref{prop:pi0}
yields as a special case 
\[ \und{\pi_0}(\cL_X^GH\und{R}) \cong \cL_X^G(\und{R}).\]

\begin{rem}
Let $R$ be a cofibrant  connective genuine
commutative $G$-ring spectrum and let $X$ be a finite simplicial
$G$-set. If we knew that the Loday construction, $\cL_X^G(R)$, were a proper
simplicial $G$-spectrum, then by an equivariant analogue of
\cite[Theorem X.2.9]{ekmm}, as used for instance in \cite[Theorem 5.2]{bghl} and \cite[Theorem 6.20]{akgh},
we would get for our  simplicial connective $G$-spectrum $\cL_X^G(R)$
and $\und{\pi}_*$ as a $G$-homology theory a spectral
sequence of the form 
\begin{equation}
\label{eq:ekmmspsq}
  E^2_{p,q} = H_p\und{\pi}_q(\cL_X^G(R)) \Rightarrow
  \und{\pi}_{p+q}|\cL_X^G(R)|\end{equation}
together with an edge homomorphism
\[ \und{\pi}_p|\cL_X^G(R)| \ra H_p\und{\pi}_0\cL_X^G(R). \]

By Proposition \ref{prop:pi0} the target of the map can be identified
with $H_p\cL_X^G(\und{\pi}_0(R))$. So we would get a  map  of
$G$-Mackey functors 
\begin{equation} \label{eq:linearization}
  \und{\pi}_k |\cL_X^G(R)| \ra \pi_k \cL_X^G(\und{\pi}_0(R))\end{equation}
that would rightly be called a linearization map.

To this end, we would need that the degeneracy maps behave well in the
sense that the map from the $n$th latching object to the $n$th
simplicial object is a cofibration for all $n$. This is true in the unstable
context if all degeneracies are $G$-cofibrations by \cite[Lemma 1.11]{mmo}.

The  degeneracies induce injective maps
  $s_i \colon X_n \ra X_{n+1}$, so the building blocks for them are isomorphisms
  $G/H \ra G/g^{-1}Hg$ and $\varnothing \ra G/H$. The first type of
  map induces a cofibration, and the second type corresponds to $S \ra
  N_H^Gi_H^*R$ that we can factor as $S \ra N_H^Gi_H^*S \ra
  N_H^Gi_H^*R$.
  As $N_H^G$ is a left Quillen functor, it preserves
  cofibrations. The restriction functors $i_H^*$ preserve cofibrations
  in the underlying category by \cite[Lemma V.2.2]{mm}. These results
  point in the right direction, but we were not
  able to find a reference that ensures properness so that the spectral sequence
  exists as in \eqref{eq:ekmmspsq}  so that a linearization map as in
  \eqref{eq:linearization} can be deduced.  
\end{rem}

\section{Spheres and suspensions} \label{sec:examples}

In the non-equivariant setting the Loday construction for the circle 
$X_\cdot=S^1$
gives Hochschild homology. We describe the equivariant Loday construction for
some circles with group action and for some unreduced suspensions
where we either flip the suspension apices with a $C_2$-action or we
fix them. We also discuss the relationship of Real topological
Hochschild homology to our equivariant Loday construction.

\subsection{Circle with rotation action}
Let $C_n$ be the cyclic group of order $n$,  $C_n = \langle \gamma
\rangle$. We let $C_n$ act on the circle $S^1_{\rot}$ by letting
$\gamma$ induce a rotation by $2\pi/n$. Then $S^1_{\rot}$ 
 has a simplicial model with non-degenerate cells being one free $0$-cell $C_n \cdot x_0
 = \{x_0, \gamma x_0, \cdots, \gamma^{n-1}x_0\}$
 and one free $1$-cell $C_n \cdot e_0$. 
 
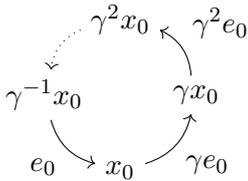
\begin{minipage}[h][][c]{0.25\linewidth}
  \begin{tikzpicture}[auto, bend right]
  \node (a) at (0,-1) {$x_0$};
  \node (b) at (1,0) {$\gamma x_0$};
  \node (c) at (0,1) {$\gamma^2 x_0$};
  \node (d) at (-1,0) {$\gamma^{-1}x_{0}$};
  \draw[->] (a) to node[swap] {$\gamma e_0$}  (b);
  \draw[->] (b) to node[swap] {$\gamma^2 e_0$} (c);
  \draw[->, dotted] (c) to  (d);
  \draw[->] (d) to node[swap] {$e_0$} (a);
\end{tikzpicture}
\end{minipage}
\begin{minipage}[h][][c]{0.7\linewidth}
We have $(S^1_{\rot})_k = \{ C_n \cdot x_k^0 , C_n \cdot x_k^1,
\cdots,  C_n \cdot x_k^k\}$,
where
\begin{align*}
  x_k^0= s_0^{k}x_0 ,
 \, x_k^i=s_0^{i-1}s_1^{k-i}e_0 \text{ for } 1 \leq i \leq k.
\end{align*}
The simplicial identities imply that
\begin{align*}
  d_j(x_k^0) = & x_{k-1}^0, \\
  d_j(x_k^i)  = & \begin{cases}
                    x_{k-1}^{i-1} & 0 \leq j \leq i-1 \\ 
                    x_{k-1}^i & i \leq j \leq k \text{ and } i \neq k
                  \end{cases} \\
  d_k(x_k^k) = & \gamma^{-1} x_{k-1}^0. 
\end{align*}
\end{minipage}

\noindent
So for a $C_n$-Tambara functor $\und{R}$ with $R := i_e^*\und{R}$, there is
\[\mathcal{L}^{C_n}_{S^1_{\rot}}(\und{R})_k = \bigsquare_{0\leq i \leq k}(C_n \otimes \und{R}) = (N_e^{C_n} R)^{\Box
  (k+1)},\]
and 
$d_i\colon (N_e^{C_n} R)^{\Box  (k+1)} \to (N_e^{C_n} R)^{\Box k}$
is 
\begin{align*}
  d_{i} & = \mathrm{id}^{i} \Box \mu \Box \mathrm{id}^{k-i} & \text{ for } 0 \leq i < k \\
  d_{k} &= (\mu \Box \mathrm{id}^{k-1}) \circ (\gamma^{-1} \Box \mathrm{id}^k) \circ \tau
\end{align*}
where $\mu\colon  (N_e^{C_n} R)^{\Box 2} \to N_e^{C_n} R$ is the multiplication and $\tau \colon
(N_e^{C_n} R)^{\Box (k+1)} \to (N_e^{C_n} R)^{\Box (k+1)}$ moves the last coordinate to
the front.

As $i_e^*\und{R}$ is an $e$-Tambara functor, it can be identified with its value on $e/e$ and that is $\und{R}(C_n/e)$.

We can identify the Loday construction with the twisted cyclic nerve
$\und{\mathrm{HC}}^{C_n}$
defined in \cite[Definition 2.20]{bghl} and its free level corresponds to a subdivision of the ordinary Loday construction. 
\begin{thm}

The $C_n$-equivariant Loday construction for  $S^1_{\rot}$ is  
\begin{equation}
\label{eq:compare-rotate-hc}
\cL^{C_n}_{S^1_{\rot}}(\und{R}) \cong \und{\mathrm{HC}}^{C_n} (N_e^{C_n} i_e^*\und{R}).
\end{equation}
\end{thm}
\begin{proof}
The claim follows by direct inspection of \cite[Definition 2.20]{bghl}. 
\end{proof}

\begin{rem}
  Note that in the Loday construction for $S^1_{\rot}$ we don't use the full structure of a $C_n$-Tambara functor, because the multiplicative norm maps are not used at all. As all simplices in $S^1_{\rot}$ are cyclically ordered, one can actually use associative Green functors instead of $C_n$-Tambara functors. This is
  the setting of \cite{bghl}. 
\end{rem}

The isomorphism \eqref{eq:compare-rotate-hc} can be generalized to the
relative case:
\begin{prop}
     Let $K \leq  C_n$ be a finite subgroup of $S^1$ and let
     $S^1_{\rot}/K$ be the circle with rotation action by $C_n$ such
     that the action by $K$ is fixed. Then the
  $C_n$-equivariant Loday construction on $S^1_{\rot}/K$ can be
  identified with the  $C_n$-twisted cyclic nerve relative to $K$ of
  \cite[Definition 3.19]{bghl}:
\begin{equation}
\label{eq:compare-rotate-hc}
\cL^{C_n}_{S^1_{\rot}/K}(\und{R}) \cong \und{\mathrm{HC}}^{C_n}_K (i_K^*\und{R}).
\end{equation}
\end{prop}
In particular, taking $K=C_n$ or $K=e$, there are isomorphisms 
\begin{align*}
\cL^{C_n}_{S^1_{\rot}/C_n}(\und{R}) & \cong \und{\mathrm{HC}}^{C_n}_{C_n} (\und{R}) \cong
\und{\mathrm{HC}}^{C_n}(\und{R}),\\
\cL^{C_n}_{S^1_{\rot}}(\und{R})  & \cong \und{\mathrm{HC}}^{C_n}_e (i_e^*\und{R}) \cong
  \und{\mathrm{HC}}^{C_n} (N^{C_n}_ei_e^*\und{R}).
\end{align*}
\begin{proof}
 Again, we choose a generator $\gamma$ so that $C_n = \langle \gamma \rangle$ and take the following simplicial model of $S^1_{\rot}/K$.
 The non-degenerate cells are one orbit of $0$-cells

\[ C_n/K \cdot x_0
 = \{Kx_0, \gamma K x_0, \cdots, \gamma^{|C_{n}/K|-1}Kx_0 = \gamma^{-1}Kx_0\}\] 
and one orbit of $1$-cells $C_n/K \cdot e_0$.

\bigskip

  \begin{minipage}[h][][c]{0.25\linewidth}
  \begin{tikzpicture}[auto, bend right]
  \node (a) at (0,-1) {$Kx_0$};
  \node (b) at (1,0) {$\gamma K x_0$};
  \node (c) at (0,1) {$\gamma^2K x_0$};
  \node (d) at (-1,0) {$\gamma^{-1}K x_{0}$};
  \draw[->] (a) to node[swap] {$\gamma K e_0$}  (b);
  \draw[->] (b) to node[swap] {$\gamma^2 K e_0$} (c);
  \draw[->, dotted] (c) to  (d);
  \draw[->] (d) to node[swap] {$K e_0$} (a);
\end{tikzpicture}
\end{minipage}
\begin{minipage}[h][][c]{0.7\linewidth}
We have $(S^1_{\rot})_k = \{ C_n/K \cdot x_k^0 , C_n/K \cdot x_k^1,
\cdots,  C_n/K \cdot x_k^k\}$,
where
\begin{align*}
  x_k^0= s_0^{k}x_0 ,
 \, x_k^i=s_0^{i-1}s_1^{k-i}e_0 \text{ for } 1 \leq i \leq k.
\end{align*}

So for a $C_n$-Tambara functor $\und{R}$ we obtain 
\[\mathcal{L}^{C_n}_{S^1_{\rot}/K}(\und{R})_k = \bigsquare_{0\leq i \leq k}(C_n/K \otimes \und{R}) =
  (N_K^{C_n} i_{K}^{*} \und{R})^{\Box
  (k+1).}\]
\end{minipage}

\bigskip \noindent 
This is precisely $\und{\mathrm{HC}}^{C_n}_K (i_K^*\und{R})_k$. 
The compatibility of this identification with the simplicial structure
maps can be seen similarly to the absolute case.
\end{proof}

If we choose the action of $\gamma^{-1}$
on $(N_e^{C_n} i_e^*\und{R})(C_n/e) = \und{R}(C_n/e)^{\otimes n}$ such that it
brings the last coordinate to the
front (see Remark \ref{rem:Weyl}), then we can identify the free-orbit level of the equivariant Loday construction with the $n$-fold subdivision (\cite[\S 1]{bhm}) of the Loday construction for the commutative ring $\und{R}(C_n/e)$: 

\begin{thm}
There is  an isomorphism
\begin{equation}
  \label{eq:compare-rotate-sd}
  \cL^{C_n}_{S^1_{\rot}}(\und{R})(C_n/e) \cong \mathrm{sd}_n\cL_{S^1} (\und{R}(C_n/e)). 
\end{equation}
\end{thm}
\begin{proof}
  Note that by the definition of the norm at the free level we obtain that 
  \[ \cL_{S^1_{\rot}}(\und{R})_k(C_n/e) = (\und{R}(C_n/e)^{\otimes
      n})^{\otimes k+1}. \] 

  We send an element
  \[ (r_{0,1} \otimes r_{0,2} \otimes \cdots \otimes r_{0,n}) \otimes \cdots \otimes (r_{k,1} \otimes \cdots \otimes r_{k,n}) \in
  \cL_{S^1_{\rot}}(\und{R})_k(C_n/e) = (\und{R}(C_n/e)^{\otimes n})^{\otimes k+1}\] to
\[ (r_{0,1} \otimes r_{1,1} \otimes\cdots \otimes r_{k,1}) \otimes (r_{0,2} \otimes \cdots \otimes
r_{k,2}) \otimes \cdots \otimes (r_{0,n} \otimes \cdots \otimes r_{k,n}) \in (\mathrm{sd}_n\cL_{S^1} \und{R}(C_n/e))_k. \] 

We have to check the compatibility of this isomorphism with the simplicial structure maps. The only non-trivial step is to compare
\begin{equation*}
  \begin{split}
    &\quad d_k((r_{0,1} \otimes r_{0,2} \otimes \cdots \otimes r_{0,n}) \otimes \cdots \otimes (r_{k,1} \otimes \cdots \otimes r_{k,n})) \\
    &= (r_{k,n}r_{0,1} \otimes r_{k,1} r_{0,2} \otimes \cdots \otimes r_{k,n-1}r_{0,n}) \otimes (r_{1,1} \otimes \cdots \otimes r_{1,n}) \otimes \cdots \otimes (r_{k-1,1} \otimes \cdots \otimes r_{k-1,n})
  \end{split}
\end{equation*}
and
\begin{equation*}
  \begin{split}
    &\quad d_k((r_{0,1} \otimes r_{1,1} \otimes\cdots \otimes r_{k,1}) \otimes (r_{0,2} \otimes \cdots \otimes
    r_{k,2}) \otimes \cdots \otimes (r_{0,n} \otimes \cdots \otimes r_{k,n})) \\
    &= (r_{k,n}r_{0,1} \otimes r_{1,1} \otimes\cdots \otimes r_{k-1,1}) \otimes (r_{k,1}r_{0,2} \otimes \cdots \otimes
r_{k-1,2}) \otimes \cdots \otimes (r_{k,n-1}r_{0,n} \otimes \cdots \otimes r_{k-1,n}).
  \end{split}
\end{equation*}
The isomorphism maps the first term to the second one. 
\end{proof}

\begin{rem}
  A similar relationship between the $H$-relative topological Hochschild homology
  $\THH_H$ for $H \leq C_n$ defined in \cite{abghlm} and the Loday
  construction for $S^1_\rot$ can be proven in the setting of 
 $C_n$-equivariant commutative ring spectra.  
\end{rem}
\subsection{Circle with reflection action}
\label{sec:circle-with-reflection}
Let $S^{\sigma}$ be the circle with reflection action. It has a $C_2$-simplicial
model with non-degenerate cells being two trivial $0$-cells and one free
$1$-cell.

\begin{minipage}[h][][c]{0.25\linewidth}
\[ \xymatrix{
    & x_0 & \\
    & & \\
 & x_1 \ar@/^5ex/[uu]^{\gamma e_0} \ar@/_5ex/[uu]_{e_0}& 
  }\] 
\end{minipage}
\begin{minipage}[h][][c]{0.7\linewidth}
We have $(S^\sigma)_k = \{ x_k^0 , C_2 \cdot x_k^1,
\cdots,  C_2 \cdot x_k^k, x_k^{k+1}\}$,
where
\begin{align*}
  x_k^0= s_0^{k}x_0 ,\, x_k^{k+1} = s_0^kx_1,\,
  x_k^i=s_0^{i-1}s_1^{k-i}e_0 \text{ for } 1 \leq i \leq k.
\end{align*}
The simplicial identities imply that
\begin{align*}
  d_j(x_k^0) = & x_{k-1}^0, \\
  d_j(x_k^{k+1}) = & x_{k-1}^k,\\
  d_j(x_k^i)  = & \begin{cases}
                    x_{k-1}^{i-1} & 0 \leq j \leq i-1, \\ 
                    x_{k-1}^i & i \leq j \leq k.
                  \end{cases}
\end{align*}
\end{minipage}

So for a $C_2$-Tambara functor $\und{R}$, there is an isomorphism of $C_2$-simplicial Tambara functors between the Loday construction and  the bar construction
\begin{equation}
  \label{eq:compare-flip}
\mathcal{L}^{C_2}_{S^{\sigma}}(\und{R}) \cong \mathrm{B}(\und{R}, N^{C_2}_ei^{*}_e\und{R}, \und{R}),
\end{equation}
where the $C_2$-Tambara structure of $\und{R}$ endows $\und{R}$ with an
$N^{C_2}_e i_e^*\und{R}$-algebra structure, coming from collapsing the free orbit to the trivial one. 

\begin{prop}
Assume that $R$ is a commutative solid ring, \ie, that the
multiplication map $\mu \colon R \otimes R \ra R$ induces an
isomorphism. If $2$ is invertible in $R=i_e^*\und{R}^c$, then $\und{R}^c$ is a
projective $N_e^{C_2}i_e^*(\und{R}^c)$-module and 
\begin{equation*}
\und{\pi}_{*}\mathcal{L}^{C_2}_{S^{\sigma}}(\und{R}^{c}) \cong \und{R}^c \Box_{N^{C_2}_ei^{*}_e(\und{R}^c)} \und{R}^c
\end{equation*}
concentrated in degree zero.
\end{prop}

\begin{proof}
  We use the explicit formula
\[  N^{C_2}_ei^{*}_e\und{R}^c =
N^{C_2}_e R =
  \begin{cases}
    \big(\mathbb{Z}\{R\} \oplus (R \otimes R)/ \mathrm{Weyl}\big)/ \mathrm{TR} & \text{ at } C_2/C_2 \\
    R \otimes R & \text{ at } C_2/e
  \end{cases}\] in order to construct a 
  splitting $$\sigma \colon  \und{R}^c \ra N^{C_2}_e R$$ of the map $N^{C_2}_e R   \to \und{R}^c$ by sending $r \in
  R$ to $r \otimes 1 \in R \otimes R$ at the $C_2/e$-level and sending $r \in R$ to
  $\left[\frac{r \otimes 1}{2}\right] \in (R \otimes R)/\mathrm{Weyl}$ at the $C_2/C_2$-level.

  We need to show that this is a morphism of  $N^{C_2}_e R$-modules
  and to this end we have to understand the multiplication in
  $N^{C_2}_e R$. We denote a generator belonging to $a \in R$ in
  $\Z\{R\}$ by $N(a)$. By Frobenius reciprocity we obtain
  \[ N(a)\cdot [x  \otimes y] = N(a) \trace(x \otimes y) =
    \trace(\res(N(a)) \cdot (x \otimes y)) = \trace((a \otimes a) \cdot (x
    \otimes y)) = [ax \otimes ay]. \]
Similarly, we get $[a \otimes b] \cdot [x \otimes y] = [ax \otimes by]
+ [bx \otimes ay]$. As the norm is multiplicative and as $N(a) =
\norm(a \otimes 1)$, we have $N(a) \cdot N(b) = N(ab)$.

At the $C_2/e$-level the map $\sigma$ is compatible with the
$N_e^{C_2}(R)$-module structure since $R$ is solid.  
At the $C_2/C_2$-level, an element
$N(r_1) \otimes r_2$ is mapped by $\id \otimes \sigma$ to $N(r_1)
\otimes [\frac{r_2\otimes 1}{2}]$ and the module action sends this
to $[\frac{r_1r_2\otimes r_1}{2}]$.
Applying first the module action, however, yields $r_1r_2r_1$ and
$\sigma$ maps this to $[\frac{r_1r_2r_1\otimes 1}{2}]$. If $R$ is a
solid commutative ring, the multiplication in $R$ identifies $
r_1r_2r_1\otimes 1$ and $r_1r_2\otimes r_1$ with each other.

Similarly, $[a\otimes b]\otimes r_2$ is mapped by $\id \otimes \sigma$ to 
$[a\otimes b]\otimes [\frac{r_2\otimes 1}{2}]$ which the module action sends to
$ [\frac{ar_2\otimes b}{2}]+  [\frac{br_2\otimes a}{2}]$ and since $R$ is a solid commutative ring this agrees with
\[ \sigma([a \otimes b]r_2) = \sigma(ar_2b + br_2a).\]
Hence  $\und{R}^c$ is a
projective $N_e^{C_2}i_e^*(\und{R}^c)$-module.
  The section $\sigma$ gives rise to a contraction of
  $B(\und{R}, N_e^{C_2}i_e^*\und{R}, \und{R})$ by sending
  $\und{R} \Box  N_e^{C_2}i_e^*\und{R}^{\Box n} \Box \und{R}$
  in simplicial degree $n$ with the map $\eta \Box \sigma \Box \id^{\Box n+1}$
  to degree $n+1$. Here $\eta \colon \und{A} \ra \und{R}$ is the unit map of
  $\und{R}$.
\end{proof}

\begin{rem}
Bousfield and Kan classified all  solid commutative  rings
\cite{bksolid}. Typical building blocks are rings of the form $\Z/n\Z$ or
subrings of the rationals $\Z[J^{-1}]$ for some set of primes $J$. 
\end{rem}

We will come back later to the $C_2$-Loday construction on $S^\sigma$, when we identify Real topological Hochschild homology with a suitable equivariant Loday construction in Theorem \ref{thm:thr}.   
\subsection{Unreduced suspension of a $G$-simplicial set}
  Let $SY_{\cdot}$ be the unreduced suspension of a finite $G$-simplicial set $Y_{\cdot}$. Using the standard
  simplicial model $\Delta^1_k = \{x_k^0, \cdots, x_k^{k+1}\}$ with $d_j(x_k^0)=x_{k-1}^0$,
  $d_j(x_k^{k+1})=x_{k-1}^{k}$ for all $0\leq j \leq k$,
  and
  \[d_j(x_k^i) = \begin{cases}
                  x_{k-1}^{i-1}, & \text{ if } 0 \leq j \leq i-1, \\
                  x_{k-1}^i, & \text{ if } i \leq j \leq k,
                \end{cases}\]
for $1 \leq i \leq k$ we get $SY_k = \{x_{k}^1 , \cdots, x_{k}^{k}\} \times Y_k \cup \{x_k^0, x_k^{k+1}\}$ and hence for any $G$-Tambara functor $\und{R}$, 
\begin{equation*}
SY_k \otimes \und{R} \cong \und{R} \Box (Y_k \otimes \und{R}) \Box \und{R}.
\end{equation*}
Keeping track of the structure maps shows that
$\mathcal{L}_{SY_{\cdot}}(\und{R})$ is isomorphic to the
diagonal of the bisimplicial Tambara functor
$\mathrm{B}(\und{R}, \mathcal{L}_{Y_{\cdot}}(\und{R}), \und{R})$.

\subsection{Unreduced suspension of a $C_2$-simplicial set} 
  Let $S^{\sigma}Y_\cdot$ be the unreduced suspension of a finite
  $C_2$-simplicial set 
  $Y_\cdot$ so that $C_2$ flips the suspension coordinate.
  The interval $[-1,1]$ with reflection action has model as a $C_2$-simplicial set 
  $[-1,1]_k = \{x_k^0, C_2\cdot x_k^1, \cdots, C_2 \cdot x_k^k, C_2\cdot x_k^{k+1}\}$ with 
  $d_j(x^0_k)=x^0_{k-1}$  and $d_j(x^{k+1}_k)=x^k_{k-1}$  for all $0\leq j \leq k$,  
  and
  \[d_j(x_k^i) = \begin{cases}
                  x_{k-1}^{i-1}, &  \text{ if } 0 \leq j \leq i-1, \\
                  x_{k-1}^i, & \text{ if } i \leq j \leq k,
                \end{cases} \, \text{ for }  1 \leq i \leq k. \]
                
We therefore get 
\begin{equation*}
S^{\sigma}Y_k = \{x_k^0\} \times Y_k \cup \{C_2\cdot x_{k}^1 , \cdots, C_2 \cdot x_{k}^{k}\} \times Y_k \cup
\{C_2 \cdot x_k^{k+1}\}
\end{equation*}
and for a $C_2$-Tambara functor $\und{R}$, 
\begin{equation*}
SY^{\sigma}_k  \otimes \und{R} = (Y_k \otimes \und{R}) \Box ((C_2
  \times Y_k) \otimes \und{R}) \Box (N^{C_2}_ei_e^*\und{R}).
\end{equation*}
Again we get that
$\cL^{C_2}_{S^{\sigma}Y_\cdot}(\und{R})$ is isomorphic to the
diagonal of the bisimplicial Tambara functor
$\mathrm{B}(\cL^{C_2}_{Y_\cdot}(\und{R}), \cL^{C_2}_{C_2 \times Y_\cdot}(\und{R}), \cL^{C_2}_{C_2}(\und{R}))$.
Note also that $S^{\sigma}Y_\cdot$ is the simplicial join  $C_2 \star Y_{\cdot}$, which is
also the homotopy pushout of $Y_\cdot \leftarrow C_2 \times Y_\cdot \to C_2$.

\subsection{Real topological Hochschild homology} \label{subsec:thr}

Hesselholt and Madsen developed Real algebraic K-theory, a variant of
algebraic K-theory that accepts as input algebras with
anti-involution \cite{hesselholt-madsen}. The corresponding Real variant of
topological Hochschild
homology, $\thr$,  was investigated in Dotto's thesis and in
\cite{dmpr} where the authors also identified $\thr(A)$ in good cases with a two-sided bar construction \cite[Prop. 2.11, Theorem 2.23]{dmpr} analogous to \eqref{eq:compare-flip}. Horev proved a similar result in the context of equivariant factorization homology \cite[Proposition 7.11]{horev}. Angelini-Knoll, Gerhardt, and Hill \cite[Definitions 4.2 and
4.5]{akgh} constructed two $O(2)$-equivariant spectra: the norm $N_{C_2}^{O(2)}(A)$
and the tensor $A \otimes_{C_2} O(2)$ using a simplicial model of $O(2)$  for a genuine commutative
$C_2$-ring spectrum $A$. They showed \cite[Propositions 4.6 and 4.9]{akgh} that
there are (zig-zag of) maps of $O(2)$-spectra
 $\thr(A) \simeq N_{C_2}^{O(2)}A$ and $N_{C_2}^{O(2)}(A) \to A \otimes_{C_2} O(2)$ such
 that the first one is a $C_2$-equivalence when $A$ is flat (\cite[Definition 3.22]{akgh}) and that the second one
 is a $C_2$-equivalence when $A$ is well-pointed (\cite[Definition 3.24]{akgh}).

We claim that there is an equivalence of simplicial
 $C_2$-spectra
\begin{equation}
\label{eq:THRequi}
 A \otimes_{C_2}  O(2)_{\bullet} \simeq \mathcal{L}^{C_{2}}_{S^{\sigma}}(A).
\end{equation}
 In fact, writing $D_{2n}$ for the dihedral group of order $2n$, so that $C_2 = D_2$,
 the $k$-simplices of $O(2)_\bullet$ are given by
 $O(2)_k = D_{4k+4}$ viewed as a $D_2$-set
 \cite[Definition 4.4]{akgh}. 
As $D_{4k+4} = \mu_{2k+2} \rtimes D_2$, we have a split short exact sequence of groups
\[ \xymatrix@1{1 \ar[r] & \mu_{2k+2} \ar[r] & D_{4k+4} \ar[r] & D_2
    \ar[r] & 1   } \]
and the induced $D_2$-action on $\mu_{2k+2}$  maps the generator
$\zeta = (1,2,\ldots,2k+2)$ to its inverse. The $D_2$-action on $D_{4k+4}$ is free and
as $D_2$-sets $D_{4k+4}/D_2 \cong \mu_{2k+2}$.
Then, 
\[ A \otimes_{D_2} D_{4k+4} \cong A \otimes \mu_{2k+2}.\]
If we choose an ordering
of the $D_2$-set $\mu_{2k+2}$ as $1 < \zeta < \zeta^2 < \ldots < \zeta^{2k+1}$, then
we always get two trivial orbits generated by $1$ and
$\zeta^{k+1}$ and $k$ free orbits generated by $\zeta, \ldots,
\zeta^k$. Hence we get that
\[ A \otimes \mu_{2k+2} \cong \mu_{2k+2} \otimes A\]
where now the tensor product of $\mu_{2k+2}$ with $A$ is the one that uses that genuine commutative $C_2$-spectra are $C_2$-commutative monoids \cite{hill-overview}.

We can identify $\mu_{2k+2}$  with the $k$-simplices of the reflection circle
$S^{\sigma}$ in Section
\ref{sec:circle-with-reflection} and this identification is compatible with the simplicial structure maps.   This shows \eqref{eq:THRequi}.

\begin{rem}
  Here, we view $D_{4k+4}$ only as a $D_2$-set. In \cite{akgh}, the group
  structure of $D_{4k+4}$ is used to set up $A \otimes_{D_2} O(2)_{\bullet}$ as a Real cyclic
  object in $C_2$-spectra. This way, the geometric realization becomes an
  $O(2)$-spectrum indexed on a $S^1$-trivial $O(2)$-universe, and $A \otimes _{D_2} O(2)$ is
  defined to be this realization after changing to a complete universe.
\end{rem}
The following result summarizes the above arguments.  
\begin{thm} \label{thm:thr} 
If $A$ is a flat
and well-pointed  $C_2$-commutative ring spectrum, then there is an equivalence of $C_2$-spectra 
\begin{equation*}
\thr(A) \simeq A \otimes_{C_2}  O(2) \simeq | \mathcal{L}^{C_{2}}_{S^{\sigma}}(A)|.
\end{equation*}
\end{thm}
\section{Relative equivariant Loday constructions} \label{sec:rel}
    In the non-equivariant context the Loday construction from
    \eqref{eq:ordinaryloday} has a relative variant: if $A$ is a commutative $k$ algebra for $k$ an arbitrary commutative ring, we can define $\cL^k_{X_\cdot}(A)$ by setting
\begin{equation} \label{eq:ordinaryrelativeloday}
 \cL^k_{X_\cdot}(A) =\{[n]\mapsto \bigotimes_{x\in X_n, k} A \} 
\end{equation}
so the tensor product over the integers is replaced by the tensor product over $k$.

Assume that $f \colon \und{R} \to \und{T}$ is a map of $G$-Tambara functors.
In the equivariant context it does \emph{not} work to replace the $\Box$-product in Definition \ref{def:equivariantloday} by the relative box product $\Box_\und{R}$. The norm terms $N_H^Gi_H^*\und{T}$ for instance don't carry an $\und{R}$-module structure in general. We propose the following definition.

\begin{defn} \label{def:relativeloday}
  Let $G$ be a finite group, $\und{R}$ and $\und{T}$ be two $G$-Tambara functors, and let  $f \colon \und{R} \to \und{T}$ be a map of Tambara functors. Then
 for any $G$-simplicial set $X_\cdot$ we define  the equivariant Loday construction of $\und{T}$ relative to $\und{R}$ on $X_\cdot$ as  
 \[\cL^{G, \und{R}}_{X_\cdot}(\und{T}) =  \cL^{G}_{X_\cdot}(\und{T}) \Box _{\cL^{G}_{X_\cdot}(\und{R}) }\und{R}.  \]
 Here, the map $\cL^{G}_{X_\cdot}(\und{R}) \to \cL^{G}_{X_\cdot}(\und{T})$ is induced by $f$
 and the map $\cL^{G}_{X_\cdot}(\und{R})\to \und{R}$ is induced by sending $X_\cdot$ to a point. The box product over a Tambara functor is defined as the usual coequalizer.
    \end{defn}
  
  \smallskip 
In \cite[Definition 1.7]{abghlm} and \cite[\S 8]{bghl} the authors define a relative norm in a similar manner.     

  In spectra, if we have a cofibrant commutative $S$-algebra $A$ and a cofibrant commutative $A$-algebra $B$, the Loday construction of $B$ on a simplicial set $X_\cdot$ is defined as 
 \[\cL_{X_\cdot}(B) = \{ [n] \mapsto \bigwedge_{x\in X_n} B\} \] 
   and the Loday construction of $B$ on $X_\cdot$ over $A$ is defined by replacing the smash products with smash products over $A$,
\[ \cL^A_{X_\cdot}(B) = \{ [n] \mapsto \bigwedge_{x\in X_n, A} B\}. \]
In \cite[\S 3]{hhlrz} we show that 
\[\cL^A_{X_\cdot}(B) \simeq \cL_{X_\cdot}(B) \wedge_{\cL_{X_\cdot}(A)} A \]
which is analogous to the definition above in the equivariant setting.

\begin{rem} \label{rem:projectivity}

  Assume that $f \colon \und{R} \ra \und{T}$ is a morphism of
  $G$-Tambara functors that turns $\und{T}$ into a projective
  $\und{R}$-module, hence it is a retract of a free 
$\und{R}$-module. As the norm functor $N_H^G$ preserves free modules
\cite[Proposition 4.0.5]{hmq} and as $N_H^G$ -- as a functor -- sends retracts
to retracts and takes values in an abelian category, $N_H^G(\und{T})$ is a
projective $N_H^G(\und{R})$-module and
in this sense, Definition \ref{def:relativeloday} is derived. 
\end{rem}

In some cases we can identify the relative Loday construction with an absolute one: 
 \begin{prop}
    For any finite group $G$ and commutative ring $R$  the map of constant Tambara functors $\und\Z^c\to \und R^c$ induced by the unit map $\Z\to R$ gives rise to  an isomorphism
\[ \cL^{G, \und\Z^c}_{X_\cdot}(\und R^c) \cong  \cL^{G }_{X_\cdot}(\und R^c) \Box \und\Z^c\] 
 for any finite simplicial $G$-set $X_\cdot$   
 \end{prop}

\begin{proof}
 As before, since $R\cong\Z\otimes R$,  by  \cite[Lemma 5.1]{lrz} we get that $\und{R}^c \cong \und{\Z}^c \Box \und{R}^c$. We then apply part (1) of  Proposition \ref{properties} to deduce $ \cL^{G}_{X_\cdot}(\und R^c) \cong  \cL^{G}_{X_\cdot}(\und R^c) \Box  \cL^{G}_{X_\cdot}(\und \Z^c)$ where the map induced by the unit map $\Z\to R$ exactly sends  $\cL^{G}_{X_\cdot}(\und \Z^c)$ to the second factor.
\end{proof}

Note that in the case when $G=C_p$ for a prime $p$ and
 $X_\cdot^{C_p} \neq \varnothing$, by Corollary \ref{Zalg} and the result
 above this in fact means that 
\begin{equation}
 \cL^{C_p, \und\Z^c}_{X_\cdot}(\und R^c) \cong  \cL^{C_p }_{X_\cdot}(\und R^c) .
 \end{equation}

 We also get a relative version of Proposition \ref{prop:trivial-constant}:

 \begin{prop}
\label{prop:relative-trivial-constant}
  Let $X_\cdot$ be a simplicial set with trivial $G$-action and let $B$ be any
  commutative $A$-algebra. Then 
\[\cL^{G,\und{A}^c}_{X_\cdot}(\und{B}^c) \cong \und{\cL^A_{X_\cdot}(B)}^c, \]
where $\cL^A_{X_\cdot}(B)$ is the nonequivariant relative Loday construction
from \eqref{eq:ordinaryrelativeloday}. 
\end{prop}
\begin{proof}
As we know by  \cite[Lemma 5.1]{lrz} that the box product of two constant Tambara functors corresponding to commutative rings is the constant Tambara functor corresponding to the tensor product of these rings, this also holds for the coequalizers $\und{B}^c \Box_{\und{A}^c} \und{B}^c \cong \und{(B \otimes_A B)}^c$ and iterations of these. So in every simplicial degree
  $n$ of 
  $\cL^{G,\und{A}^c}_{X_\cdot}(\und{B}^c)$, we have $\Box_{x\in X_n, \und{A}^c} \und{B}^c \cong \und{(\otimes_{x\in X_n, A} B)}^c$ and these isomorphisms are compatible with the simplicial structure maps.
\end{proof}


\begin{bibdiv}
\begin{biblist}


\bib{anderson}{incollection}{
    AUTHOR = {Anderson, D. W.},
     TITLE = {Chain functors and homology theories},
 BOOKTITLE = {Symposium on {A}lgebraic {T}opology ({B}attelle {S}eattle
              {R}es. {C}enter, {S}eattle, {W}ash., 1971)},
    SERIES = {Lecture Notes in Math., Vol. 249},
     PAGES = {1--12},
 PUBLISHER = {Springer, Berlin},
      YEAR = {1971},
}


\bib{akgh}{misc}{
    author={Angelini-Knoll, Gabriel},
    author={Gerhardt, Teena},
    author={Hill, Michael}, 
    title={Real topological Hochschild homology via the norm and Real
      Witt vectors}, 
   note={preprint:  arXiv:2111.06970}, 
  }


 
\bib{abghlm}{article}{
  AUTHOR = {Angeltveit, Vigleik},
  AUTHOR = {Blumberg, Andrew J.},
  AUTHOR = {Gerhardt, Teena},
  AUTHOR = {Hill, Michael A.},
  AUTHOR = {Lawson, Tyler},
  AUTHOR = {Mandell, Michael A.},
     TITLE = {Topological cyclic homology via the norm},
   JOURNAL = {Doc. Math.},
      VOLUME = {23},
      YEAR = {2018},
     PAGES = {2101--2163},
} 


\bib{bghl}{article}{
  AUTHOR = {Blumberg, Andrew J.},
  author = {Gerhardt, Teena},
  author={Hill, Michael A.}, 
  author= {Lawson, Tyler},
     TITLE = {The {W}itt vectors for {G}reen functors},
   JOURNAL = {J. Algebra},
      VOLUME = {537},
      YEAR = {2019},
     PAGES = {197--244},
}


\bib{bhm}{article}{
  AUTHOR = {B\"{o}kstedt, M.},
  author={Hsiang, W. C.},
  author={Madsen, I.},
     TITLE = {The cyclotomic trace and algebraic {$K$}-theory of spaces},
   JOURNAL = {Invent. Math.},
    VOLUME = {111},
      YEAR = {1993},
    NUMBER = {3},
     PAGES = {465--539},
   }


\bib{bksolid}{article}{
  AUTHOR = {Bousfield, Aldridge. K.},
  author = {Kan, Dan M.},
     TITLE = {The core of a ring},
   JOURNAL = {J. Pure Appl. Algebra},
      VOLUME = {2},
      YEAR = {1972},
      PAGES = {73--81},
      }

\bib{brunfree}{article}{
    AUTHOR = {Brun, Morten},
     TITLE = {Witt vectors and {T}ambara functors},
   JOURNAL = {Adv. Math.},
      VOLUME = {193},
      YEAR = {2005},
    NUMBER = {2},
    PAGES = {233--256},
  }


\bib{brunpi0}{article}{
    AUTHOR = {Brun, Morten},
     TITLE = {Witt vectors and equivariant ring spectra applied to
              cobordism},
   JOURNAL = {Proc. Lond. Math. Soc. (3)},
  FJOURNAL = {Proceedings of the London Mathematical Society. Third Series},
    VOLUME = {94},
      YEAR = {2007},
    NUMBER = {2},
     PAGES = {351--385},
  }

  \bib{bds}{misc}{
    author={Brun, Morten},
    author={Dundas, Bj{\o}rn Ian}, 
    author={Stolz, Martin},
    title={Equivariant Structure on Smash Powers}, 
    note={preprint: arXiv:1604.05939},
  }

\bib{dmpr}{article}{
  AUTHOR = {Dotto, Emanuele},
  author={Moi, Kristian},
  author={Patchkoria, Irakli},
  author={Reeh, Sune Precht},
     TITLE = {Real topological {H}ochschild homology},
   JOURNAL = {J. Eur. Math. Soc. (JEMS)},
      VOLUME = {23},
      YEAR = {2021},
    NUMBER = {1},
     PAGES = {63--152},
}
  
  
\bib{ekmm}{book}{
TITLE={Rings, modules, and algebras in stable homotopy theory},
AUTHOR={Elmendorf, A. D.},
AUTHOR={Kriz, I.},
AUTHOR={Mandell, M. A.},
AUTHOR={May, J.P.},
DATE={1997},
PUBLISHER={American Mathematical Society, Mathematical Surveys and Monographs 47}
ADDRESS={Providence, RI}
}


\bib{hhlrz}{article}{
  AUTHOR = {Halliwell, Gemma},
  author={H\"{o}ning, Eva},
  author= {Lindenstrauss, Ayelet}, 
  author={Richter, Birgit},
  author={Zakharevich, Inna},
     TITLE = {Relative {L}oday constructions and applications to higher
              {$\mathsf{THH}$}-calculations},
   JOURNAL = {Topology Appl.},
    VOLUME = {235},
      YEAR = {2018},
     PAGES = {523--545},
}


\bib{hesselholt-madsen}{book}{
  author={Hesselholt, Lars},
  author={Madsen, Ib},
  title={Real algebraic K-theory},
  note={preliminary version, available at \url{https://web.math.ku.dk/~larsh/papers/s05/}},
}



\bib{hill-overview}{incollection}{
    AUTHOR = {Hill, Michael A.},
     TITLE = {Equivariant stable homotopy theory},
 BOOKTITLE = {Handbook of homotopy theory},
    SERIES = {CRC Press/Chapman Hall Handb. Math. Ser.},
     PAGES = {699--756},
 PUBLISHER = {CRC Press, Boca Raton, FL},
      YEAR = {2020},
}

\bib{hh}{misc}{
  AUTHOR={Hill, Michael A.},
  AUTHOR = {Hopkins, Michael J.},
   TITLE = {Equivariant Symmetric Monoidal Structures}, 
      NOTE={preprint, arXiv:1610.03114v1}, 
    }


\bib{hhr}{article}{
  AUTHOR = {Hill, Michael A.},
  author= {Hopkins, Michael J.},
  author={Ravenel, Doug C.},
     TITLE = {On the nonexistence of elements of {K}ervaire invariant one},
   JOURNAL = {Ann. of Math. (2)},
    VOLUME = {184},
      YEAR = {2016},
    NUMBER = {1},
     PAGES = {1--262},    
}

     
  \bib{hm}{article}{
    AUTHOR={Hill, Michael A.},
    AUTHOR={Mazur, Kristen},
TITLE={An equivariant tensor product on Mackey functors}, 
JOURNAL={J. Pure Appl. Algebra},
YEAR={2019},
NUMBER = {12},
PAGES={5310--5345},
}

 \bib{hmq}{article}{
   AUTHOR={Hill, Michael A.},
   AUTHOR={Mehrle, David},
   AUTHOR={Quigley, James D.},
TITLE={Free Incomplete Tambara Functors are Almost Never Flat}, 
JOURNAL={International Mathematics Research Notices},
YEAR={2023},
NUMBER = {5},
PAGES={4225--4291},
}



\bib{horev}{misc}{
  AUTHOR={Horev, Asaf}, 
TITLE={Genuine equivariant factorization homology}, 
NOTE={preprint, arXiv:1910.07226}, 
}

\bib{hoyer}{misc}{
AUTHOR={Hoyer, Rolf},
TITLE={Two topics in stable homotopy theory},
YEAR={2014},
NOTE={Dissertation, the University of Chicago}, }


    \bib{lrz}{article}{
      AUTHOR = {Lindenstrauss, Ayelet},
      AUTHOR = {Richter, Birgit},
      AUTHOR = {Zou, Foling},
      TITLE ={Examples of \'etale extensions of Green functors},
      YEAR={2024},
      JOURNAL={Proceedings of the AMS, Ser. B},
      NUMBER = {11},
      PAGES={287--303},
    }


    \bib{mm}{article}{
      AUTHOR = {Mandell, Michael A.},
      AUTHOR = {May, J. Peter},
     TITLE = {Equivariant orthogonal spectra and {$S$}-modules},
   JOURNAL = {Mem. Amer. Math. Soc.},
      VOLUME = {159},
      YEAR = {2002},
    NUMBER = {755},
     PAGES = {x+108},
   }


   \bib{mmo}{book}{
AUTHOR = {May, J. Peter},
      AUTHOR = {Merling, Mona},
      AUTHOR = {Osorno, Ang\'elica},
      TITLE ={Equivariant infinite loop space theory, the space level story},
      series={Mem. Amer. Math. Soc.},
      volume={305},
      year = {2025},
      number={1540},
      pages={v+136 pp.}
  }
   
    
    \bib{mazur}{misc}{
      AUTHOR={Mazur, Kristen},
      TITLE={On the structure of Mackey functors and Tambara functors}, 
      YEAR = {2013}, 
      NOTE={Dissertation, University of Virginia}, 
}


\bib{pirashvili}{article}{
    AUTHOR = {Pirashvili, Teimuraz},
     TITLE = {Hodge decomposition for higher order {H}ochschild homology},
   JOURNAL = {Ann. Sci. \'{E}cole Norm. Sup. (4)},
    VOLUME = {33},
      YEAR = {2000},
    NUMBER = {2},
     PAGES = {151--179},
}

    
\bib{cats}{book}{
    AUTHOR = {Richter, Birgit},
     TITLE = {From categories to homotopy theory},
    SERIES = {Cambridge Studies in Advanced Mathematics},
    VOLUME = {188},
 PUBLISHER = {Cambridge University Press, Cambridge},
      YEAR = {2020},
     PAGES = {x+390},
}


\bib{strickland}{misc}{
 AUTHOR = {Strickland, Neil},
  TITLE = {Tambara functors}, 
     NOTE={preprint, arXiv:1205.2516}, 
   }

   
\bib{tambara}{article}{
    AUTHOR = {Tambara, Daisuke},
     TITLE = {On multiplicative transfer},
   JOURNAL = {Comm. Algebra},
    VOLUME = {21},
      YEAR = {1993},
    NUMBER = {4},
     PAGES = {1393--1420},
 } 

 
 \bib{tw}{article}{
   AUTHOR = {Th\'{e}venaz, Jacques},
   AUTHOR = {Webb, Peter},
     TITLE = {The structure of {M}ackey functors},
   JOURNAL = {Trans. Amer. Math. Soc.},
  FJOURNAL = {Transactions of the American Mathematical Society},
    VOLUME = {347},
      YEAR = {1995},
    NUMBER = {6},
     PAGES = {1865--1961},
 }

 \bib{ullman}{misc}{
author={Ullman, John}, 
title={Tambara Functors and Commutative Ring Spectra}, 
note={preprint: arXiv:1304.4912}, 
 }

\end{biblist}
\end{bibdiv}
\end{document}